\documentclass[10pt,a4paper]{article}
\usepackage{hyperref}
\usepackage[utf8]{inputenc}
\usepackage{amsmath}
\usepackage{amsfonts}
\usepackage{amssymb}
\usepackage{makeidx}
\usepackage{graphicx}
\usepackage{hyperref}
\usepackage{amssymb,dsfont,amsthm,amsmath,makeidx,verbatim,latexsym,amsfonts,mathrsfs}
\usepackage{cancel}
\usepackage{graphicx}
\usepackage[none]{hyphenat}
\usepackage[english]{babel}
\addto{\captionsenglish}{%

}
\DeclareMathAlphabet{\mathpzc}{OT1}{pzc}{m}{it}

\theoremstyle{plain}
\newtheorem{theorem}{\textbf{Theorem}}[section]

\newtheorem{definition}[theorem]{\textbf{Definition}}

\newtheorem{thm}{Theorem}[section]
\newtheorem{prop}[thm]{Proposition}
\newtheorem{lem}[thm]{Lemma}

\newtheorem{rem}[thm]{Remark}

\numberwithin{equation}{section}

\title{\textbf{A Note on Some Classes of Function Algebras}}
\author{M. E. Egwe$^1$ and R.F. Yusuf$^2$\\ Department of Mathematics, University of Ibadan, Ibadan.\\  $^1$\emph{murphy.egwe@ui.edu.ng},\;$^2$\emph{funke.yusuff@gmail.com}}
\begin{document}
\maketitle

\large	
\begin{abstract}
In this note, an attempt is made at highlighting a construct of the algebras of holomorphic functions on the unit disc,  bi-invariant functions,  Lie pseudo groups and the Colombeau's algebra (of generalized functions). We establish some new results on these algebras and also give new proofs to some existing ones.
\end{abstract}

 \baselineskip 10pt	
\section{Introduction}
	The theory of holomorphic functions on the unit disc and its ubiquity abound in literature (See \cite{Zhu04},\cite{Zhu07}). Here, we shall concern ourselves with introducing some of the basic definitions needed in the succeeding developments of some function algebras and proving several essential algebraic results.
\begin{definition}\cite{Kaniuth}
\normalfont
 	A linear space $\mathcal{A}$ over $\mathbb{C}$ is said to be an algebra if it is equipped with a binary operation, referred to as multiplication and denoted by juxtaposition, from $\mathcal{A}\times\mathcal{A}\to\mathcal{A}$ such that
	\begin{gather}
	f(gh)=(fg)h,\\
	f(g+h)=fg+fh; \ (g+h)f=gf+hf,\\
	a(fg)=(af)g=f(ag).	
	\end{gather}
	for all $f,g,h\in\mathcal{A}$ and $a\in\mathbb{C}$.
\end{definition}	
 $\mathcal{A}$ is called a commutative algebra if $\mathcal{A}$ is an algebra and
	\begin{equation}
	fg=gf
	\end{equation}
	for all $f,g\in\mathcal{A}$, whereas $\mathcal{A}$ is an algebra with identity if $\mathcal{A}$ is an algebra and there exists some element $e\in\mathcal{A}$ such that
	\begin{equation}
	ef=fe=f
	\end{equation}
	for all $f\in\mathcal{A}$.
\begin{definition}
\normalfont
	A normed linear space$\left(\mathcal{A},\left\|\cdot\right\|\right)$ over $\mathbb{C}$ is said to be normed algebra if $\mathcal{A}$ is an algebra and
	$$\left\|fg\right\|\leq\left\|f\right\|\left\|g\right\|$$ for all $f,g\in\mathcal{A}$.
\end{definition}	
	A normed algebra $\mathcal{A}$ is said to be a Banach algebra if the normed linear space $\left(\mathcal{A},\left\|\cdot\right\|\right)$ is a Banach space.
\begin{definition}
\normalfont
	Let $\mathcal{A}$ be a complex algebra. An involution on $\mathcal{A}$ is a mapping	$* : f\to f^{*}$ from $\mathcal{A}$ into $\mathcal{A}$ satisfying the following conditions.
	\begin{gather}
		\left(f+g\right)^{*} = f^{*}+g^{*},\\
	 	\left(\lambda f\right)^{*} = \overline{\lambda} f^{*},\\
		\left(fg\right)^{*} = g^{*}f^{*},\\
		\left(f^{*}\right) = f,
	\end{gather}
	for all $f,g\in\mathcal{A}$ and $\lambda\in\mathbb{C}$. $\mathcal{A}$ is then called a $*-$algebra or an algebra with
	involution.
	
	A $\mathit{C}^{*}-$algebra is a Banach algebra $\mathcal{A}$ with involution in which for all $f\in\mathcal{A}$,
	\begin{equation}
	\left\|f^{*}f\right\|=\left\|f\right\|^{2}.
	\end{equation}
\end{definition}	
	A homomorphism between $\mathit{C}^{*}-$algebras $\mathcal{A}$ and $\mathcal{B}$ is a linear map $\varphi:\mathcal{A}\to\mathcal{B}$ that satisfies $\varphi(fg)=\varphi(f)\varphi(f)$ and $\varphi(f^{*})=\varphi(g)^{*}$ for all $f\in\mathcal{A}$ and $g\in\mathcal{B}$ \cite{Murphy}.
An isomorphism between two $\mathit{C}^{*}-$algebras is an invertible homomorphism.
	
	The dual $\mathcal{A}^{\prime}$ of a $*-$algebra $\mathcal{A}$ has a canonical involution given by $f^{*}(x)=\overline{f(x^{-1})}$ for all $f\in\mathcal{A}$.
\begin{definition}
\normalfont
	A function algebra on a compact  Hausdorff space $X$ is a  commutative Banach algebras $\mathcal{A}$ over $\mathbb{C}$ which satisfies the  following conditions:
	\begin{enumerate}
		\item  The elements of $\mathcal{A}$ are continuous complex-valued functions defined on $X$,  i.e. $\mathcal{A}\subset\mathit{C(X)}$;
		\item $\mathcal{A}$ contains all constant functions on $X$;
		\item The operations on $\mathcal{A}$ are the pointwise additions and multiplication;
		\item $\mathcal{A}$ is closed with respect to the uniform norm in $\mathit{C(X)}$
		\begin{equation*}
		\Vert f\Vert =\underset{x\in X}{\mathrm{sup}}|f(x)|, \ f\in \mathcal{A};
		\end{equation*}
		\item $\mathcal{A}$ separates the  points of $X$.
	\end{enumerate}
	Function algebras are commutative Banach algebras over $\mathbb{C}$ with unit (the constant  function $1$ on $X$) and this fact plays a crucial role in their studies \cite{Larsen}.
\end{definition}	
	Let $X$ be a compact Hausdorff space and $\mathcal{A}$ a non-empty subset of $\mathit{C(X)}$. For each $x\in X$ the evaluation map at $x$, denoted by $\phi_x$, is defined  by
	\begin{equation}
	\phi_{x}(f)=f(x) \ \mathrm{for} f\in \mathcal{A}
	\end{equation}
	We observe that if $\mathcal{A}$ is a subspace, then $\phi_{x}:\mathcal{A}\to \mathbb{C}$ is a linear map, and  if $\mathcal{A}$ is subalgebra, then $\phi_{x}$ is a homomorphism. If $\mathcal{A}$ contains $1$ then $\phi_{x}(1)=1$ and hence $\phi_{x} \neq 0$.

Unitization of normed algebras is well-known \cite{Kaniuth}.	
\begin{definition}
\normalfont
	Let $\mathcal{A}$ be a normed algebra. A left (right) approximate	identity for $\mathcal{A}$ is a net $\left(e_{\alpha}\right)_{\alpha}$ in $\mathcal{A}$ such that $e_{\alpha}f\to f$  $(fe_{\alpha} \to f)$ for each $f\in\mathcal{A}$. An approximate identity for $\mathcal{A}$ is a net $(e_{\alpha})$ which is both a left and a right approximate identity. A (left or right) approximate identity $\left(e_{\alpha}\right)_{\alpha}$ is bounded by $M>0$ if $\left\|e_{\alpha}\right\|\leq M$ for all $\alpha$.
A has left (right) approximate units if, for each $f\in\mathcal{A}$ and $\epsilon>0$, there exists $u\in\mathcal{A}$ such that $\left\|f-uf\right\|\leq\epsilon$ $(\left\|f-fu\right\|\leq\epsilon)$, and $\mathcal{A}$ has
	an approximate unit if, for each $f\in\mathcal{A}$ and $\epsilon>0$, there exists $u\in\mathcal{A}$ such
	that $\left\|f-uf\right\|\leq\epsilon$ and $\left\|f-fu\right\|\leq\epsilon$. A has a (left, right) approximate unit bounded by $M>0$, if the elements $u$ can be chosen such that $\left\|u\right\|\leq M$.
\end{definition}	
Suppose that $\mathcal{A}$ is an approximately unital Banach algebra. We	define the unitization of $\mathcal{A}$ by considering the canonical ‘left regular representation’	$\lambda: \mathcal{A}\to B(\mathcal{A})$, and identifying $\mathcal{A} + \mathbb{C}_1$ with the span of $\lambda(\mathcal{A})+\mathbb{C}I_{\mathcal{A}}$, which is a unital Banach subalgebra of $B(\mathcal{A})$. Thus if $f\in\mathcal{A}$ and $\mu\in\mathbb{C}$ then
	\begin{equation}
	\left\|f+\mu 1\right\|=\underset{\Vert c\Vert\leq 1}{\text{sup}}\left\|fc+\mu c\right\|,
	\end{equation}
	for all $c\in\mathcal{A}$.
	
	It is occasionally useful that there are some other equivalent expressions for	the quantity above. For example, if $\left(e_{\alpha}\right)_{\alpha}$ is an approximate identity for $\mathcal{A}$ then
	\begin{equation}
	\left\|f+\mu_1\right\|=\underset{\alpha}\lim\left\|fe_{\alpha}+\mu e_{\alpha}\right\|=\underset{\alpha}{\text{sup}}\left\|f e_{\alpha}+\mu e_{\alpha}\right\|.
	\end{equation}	
Given a real algebra $\mathcal{A}$, the complexification $\mathcal{B}$ of $\mathcal{A}$ is the set $\mathcal{A} \times\mathcal{A}$ with the operations of addition, multiplication, and scalar multiplication defined by
	\begin{gather*}
	\left(f_{1},g_{1}\right)+\left(f_{2},g_{2}\right):=\left(f_{1}+f_{2},g_{1}+g_{2}\right),\\
	\left(\alpha+i\beta\right)\left(f,g\right):=\left(\alpha f-\beta g,\alpha g+\beta f\right),\\
	\left(f_{1},g_{1}\right)\left(f_{2},g_{2}\right):=\left(f_{1}f_{2}-g_{1}g_{2},f_{1}g_{2}+g_{1}f_{2}\right),
	\end{gather*}
	for all $f_{1},f_{2},g_{1},g_{2}\in\mathcal{A}$ and $\alpha,\beta\in\mathbb{R}$.
\begin{definition}
\normalfont
	An operator $\varphi : \mathcal{H}_{1} \to\mathcal{H}_{2}$ between two Hilbert space is simply a linear map (i.e., $\varphi\left(\lambda f + \mu g\right) = \lambda\varphi (f) + \mu\varphi(g)$ for all $\lambda,\mu\in\mathbb{C}$ and $f,g\in\mathcal{H}_{1}$.
	
	An operator on a Hilbert space $\mathcal{H}$ is bounded if and only if it is continuous in the sense that $f_{n}\to f$ implies $\varphi f_{n} \to\varphi f$ for all convergent sequences $\left(f_{n}\right)$ in $\mathcal{H}$.
\end{definition}	
	Let $\varphi:\mathcal{H}_{1}\to\mathcal{H}_{2}$ be an operator. We define $\left\|\varphi\right\|\in\mathbb{R}^{+}\cup\left\{\infty\right\}$ by
	\begin{equation}
	\left\|\varphi\right\|:=\underset{\left\|f\right\|_{\mathcal{H}_{1}}=1}{\text{sup}}\left\|\varphi(f)\right\|_{\mathcal{H}_{2}},
	\end{equation}
	for all $f\in\mathcal{H}_{1}$, where $\left\|f\right\|_{\mathcal{H}_{1}}=\left|\left\langle f,f\right\rangle\right|_{\mathcal{H}_{1}},$ etc. We say that $\varphi$ is bounded when $\left\|\varphi\right\|<\infty$, in which case the number $\left\|f\right\|$ is called the norm of $\varphi$.
\subsection{Some Classes of Operator Algebras}
	\begin{enumerate}
\item [(1)] Operator Algebra Valued Continuous Functions: Let $X$ be a compact space, and $\mathcal{A}$ an operator algebra. Then the operator space $\mathcal{C}(X;\mathcal{A})$ is an operator algebra when equipped with the product defined by $(fg)(t)=f(t)g(t)$. For example, if $\mathcal{A}$ is a subalgebra of $\mathcal{B(H)}$ then $\mathcal{C}(X;\mathcal{A})$ is a subalgebra of the $\mathcal{C}^{*}-$algebra $\mathcal{C}(X;\mathcal{B(H)})$. If $\mathcal{A}$ is unital, then $\mathcal{C}(X;\mathcal{A})$ is also unital, the identity being the constant function equal to the identity of $\mathcal{A}$.
\item [(2)] Uniform Algebras: A (concrete) uniform algebra is a unital-subalgebra of $\mathcal{C}(X)$, for some compact space $X$. Here, we will consider any uniform algebra as endowed with its minimal operator space structure. Then an (abstract) uniform algebra is a unital operator algebra which is completely isometrically isomorphic to a concrete uniform algebra. In this way, we regard uniform algebras as a subclass of the operator algebras.
		
More generally, we will use the term \emph{function algebra} for an operator algebra $\mathcal{A}$ for which there exists a compact space $X$ and a completely isometric homomorphism $\pi:\mathcal{A}\to\mathcal{C}(X)$. Any function algebra is a minimal operator space.
\item [(3)] Disc Algebras: This fundamental example of a uniform algebra has two equivalent definitions. Let us denote by $\mathbb{B}$ and $\mathbb{T}$ the open unit disc of $\mathbb{C}$ and the unitary complex group $\mathbb{T}=\left\{z\in\mathbb{C}:|z|=1\right\}$ respectively. Then the disc algebra $\mathcal{A(\mathbb{B})}$ is the subalgebra of $\mathcal{C}(\overline{\mathbb{B}})$ consisting of all continuous functions $F:\overline{\mathbb{B}}\to\mathbb{C}$, whose restriction to $\mathbb{B}$ is holomorphic. By the maximal modulus theorem, the restriction of functions in $\mathcal{A(\mathbb{B})}$ to the boundary $\mathbb{T}$ is an isometry. Hence we may alternatively regard $\mathcal{A(\mathbb{B})}\subset\mathcal{C}(\mathbb{T})$ as a uniform algebra acting on $\mathbb{T}$. In that representation, $\mathcal{A(\mathbb{B})}$ consists of all elements of $\mathcal{C}(\mathbb{T})$ whose harmonic extension to $\mathbb{B}$ given by the Poisson integral is holomorphic. Equivalently, given any $f\in\mathcal{C}(\mathbb{T})$, we associate Fourier coefficients
		\begin{equation}
		\widehat{f}(k)=\int_{\mathbb{T}}f(z)z^{-k}d\mu(z), \ \ k\in\mathbb{Z}.
		\end{equation}
Then, $\mathcal{A(\mathbb{B})}\subset\mathit{C(\mathbb{T})}$ is the closed subalgebra of all $f\in\mathcal{C}(\mathbb{T})$ such that $f(k)=0$ for every $k<0$. The (holomorphic) polynomials form a dense subalgebra of $\mathcal{A(\mathbb{B})}$.
\item [(4)] Adjoint Algebra: The adjoint operator space $\mathcal{A}^{*}$ is an operator algebra, with product $a^{*}b^{*}=(ba)^{*}$, for $a,b\in\mathcal{A}$. Indeed, if $\mathcal{A}$ is a subalgebra of a $\mathcal{C}^{*}-$algebra $\mathcal{B}$, then $\mathcal{A}^{*}$ may be identified with the subalgebra ${a^{*}:a\in\mathcal{A}}$ of $\mathcal{B}$. Note that if $\mathcal{A}$ has a continuous approximate identity (cai) $(e_{t})_{t}$, then $(e^{*}_{t})_{t}$ is a cai for $\mathcal{A}^{*}$.
\item [(5)] Multiplier Operator Algebras \cite{Blecher}: A left (right) multiplier on $\mathcal{A}$ is a linear mapping $T : \mathcal{A}\to\mathcal{A}$ such that $$\mathit{T}(xy)=\mathit{T}(x)y \ (= x\mathit{T}(y))$$ for all $x,y\in\mathcal{A}$.
		
$\mathit{T}$ is called a two-sided multiplier (or simply, a multiplier) on $\mathcal{A}$ if it is a left and a right multiplier.
		
Let $\mathcal{A}$ be a $\mathcal{C}^{*}-$algebra. The multiplier algebra of $\mathcal{A}$, denoted by $\mathcal{M(A)}$, is the universal $\mathcal{C}^{*}-$algebra with the property that $\mathcal{M(A)}$ contains $\mathcal{A}$ as an essential ideal and for any $\mathcal{C}^{*}-$algebra $\mathcal{B}$ containing $\mathcal{A}$ as an essential ideal there exists a unique  $^{*}-$homomorphism $\pi:\mathcal{B}\to\mathcal{M(A)}$, that is, the identity on $\mathcal{A}$.
		
If $\mathcal{A}$ is a unital $\mathcal{C}^{*}-$algebra. Then $\mathcal{M(A)}$ is unital and $\mathcal{M(A)}=\mathcal{A}$.
	\end{enumerate}
	
\subsection{Bi-invariant Function}
	Let $G$ be a locally compact group and $K$ a compact subgroup. We assume the normalization $\int_{K}d\mu_{k}(k)=1$. The projection
	\begin{equation*}\label{key}
		P:G\to K\setminus G/K
	\end{equation*}
	defined by $p(g)=KgK$ identifies $C_{c}\left( K\setminus G/K\right)$ with
	\begin{equation*}
		\left\{\phi\in C_{c}\left(K\setminus G/K\right):\phi(k_{1}gk_{2})=\phi(g) \ \mathrm{for} \ k_{i}\in K, \ g\in G\right\}.
	\end{equation*}
Now, let us identify $C\left(K\setminus G/K\right), \ \mathcal{C}_{\infty}\left(K\setminus G/K\right)$, and $\mathit{L}^{p}\left(K\setminus G/K\right)$ with the bi-K-invariant functions of the same category on $G$. In the sequel, $\mathcal{C}_{\infty}$ means the commutative Banach algebra of continuous functions that vanish at infinity with pointwise multiplication and sup norm. In other words, $\mathcal{C}_{\infty}(G)$ consists of the continuous functions $f:G\to\mathbb{C}$ such that, if $\epsilon>0$ then there exists a compact set $C\subset G$ such that $\left|\phi(x)\right|<\epsilon$ for all $x\in G\setminus C$.
		
	We now have the projections $C(G)\to C\left(K\setminus G/K\right)$, $C_{c}(G)\to C_{c}\left(K\setminus G/K\right)$, $\mathcal{C}_{\infty}(G)\to\mathcal{C}_{\infty}\left(K\setminus G/K\right)$ and $\mathit{L}^{p}(G)\to\mathit{L}^{p}\left(K\setminus G/K\right)$, all denote $\phi\to\phi^{\sharp}\in C_{c}(G)^{\sharp}$ given by
	\begin{equation}
		\phi^{\sharp}(g)=\int_{K}\int_{K}\phi(k_{1}gk_{2})d\mu_{k}(k_{1})(k_{2}).
	\end{equation}
	Some immediate properties are (see \cite{Loomis}):
	\begin{itemize}
		\item If $\phi_{1}\in C\left(K\setminus G/K\right)$ and $\phi_{2}\in C(G)$ the $(\phi_{1}\phi_{2})^{\sharp}=\phi_{1}\phi_{2}^{\sharp}$ for the pointwise multiplication, and
		\item If $\phi_{1}\in C_{c}\left(K\setminus G/K\right)$ and $\phi_{2}\in C_{c}(G)$ then $\left(\phi_{2}*\phi_{1}\right)^{\sharp}=\phi_{2}^{\sharp}* \phi_{1}$ and $\left(\phi_{1}*\phi_{2}\right)^{\sharp}=\phi_{1}^{\sharp}*\phi_{2}$ for the convolution product. In particular, $C_{c}\left(K\setminus G/K\right)$ is a subalgebra of the convolution algebra $C_{c}(G)$ and $\mathit{L}^{1}\left(K\setminus G/K\right)$ is a subalgebra of $\mathit{L}^{1}(G)$.
	\end{itemize}
	
	We say that $(G,K)$ is Gelfand pair if the convolution algebra $\mathit{L}^{1}\left(K\setminus G/K\right)$ is commutative. If $(G,K)$ is a Gelfand pair, then $G/K$ is a commutative space relative to $G$, and we also say that $(G,K)$ is a commutative pair. Since $C_{c}\left(K\setminus G/K\right)$ is dense in $\mathit{L}^{1}\left(K\setminus G/K\right)$ it is equivalent to require that $C_{c}\left(K\setminus G/K\right)$ be commutative.
\subsection{Holomorphic Functions}
As usual, we write $\mathbb{N}$ for the set of natural numbers, $\mathbb{R}$ for the field of real numbers and $\mathbb{C}$ for the field of complex numbers. For any positive integer $n\in\mathbb{N}$ the set $\mathbb{C}^{n}$ is equipped with the usual vector space structure and for any $z = \left(z_{1}, \cdots , z_{n}\right)\in\mathbb{C}^{n}$ the norm of $z$ is given by the popular Euclidean norm $$|z| = \left(|z_{1}|^{2} + \cdots + |z_{n}|^{2}\right)^{1/2}.$$ We define an isomorphism of $\mathbb{R}-$vector spaces between $\mathbb{C}^{n}$ and $\mathbb{R}^{2n}$ by setting	 $z_{j} = x_{j} + iy_{j}$ for any $z = \left(z_{1}, \cdots , z_{n}\right)\in\mathbb{C}^{n}$ and $j = 1,\cdots, n$.
	
	The holomorphic and anti-holomorphic differential operations are given by
	\begin{equation}
	\begin{cases}
	&\frac{\partial}{\partial_{j}}=\frac{1}{2}\left(\frac{\partial}{\partial x_{j}}+\frac{1}{i}\frac{\partial}{\partial y_{j}}\right)=\frac{1}{2}\left(\frac{\partial}{\partial x_{j}}-i\frac{\partial}{\partial y_{j}}\right), \ j=1,\cdots,n\\
	&
	\frac{\partial}{\partial\overline{z}_{j}}=\frac{1}{2}\left(\frac{\partial}{\partial x_{j}}-\frac{1}{i}\frac{\partial}{\partial y_{j}}\right)=\frac{1}{2}\left(\frac{\partial}{\partial x_{j}}+i\frac{\partial}{\partial y_{j}}\right), \ j=1,\cdots,n.
	\end{cases}
	\end{equation}
	
	If $\alpha=\left(\alpha_{1},\cdots,\alpha_{n}\right)\in\mathbb{N}^{n}$ and $\beta=\left(\beta_{1},\cdots,\beta_{n}\right)\in\mathbb{N}^{n}$ are multi-indices and $x=\left(x_{1},\cdots,x_{n}\right)$ is a  point in $\mathbb{R}^{n}$ then we set $$\left|\alpha\right|=\alpha_{1}+\cdots+\alpha_{n}, \ \alpha!=\alpha_{1}!\cdots\alpha_{n}!, \ x^{\alpha}=x_{1}^{\alpha_{1}}\cdots x_{n}^{\alpha_{n}},$$
	\begin{equation}
	D^{\alpha}=\frac{\partial^{|\alpha|}}{\partial x_{1}^{\alpha_{1}}\cdots \partial x_{n}^{\alpha_{n}}}  \ \text{and} \ D^{\alpha\overline{\beta}}=\frac{\partial^{|\alpha|+|\beta|}}{\partial z_{1}^{\alpha_{1}}\cdots\partial z_{n}^{\alpha_{n}}\partial\overline{z}_{1}^{\beta_{1}}\cdots\partial\overline{z}_{n}^{\beta_{n}}}.
	\end{equation}
	
	If $\mathbb{B}$ is an open ball in $\mathbb{R}^{n}$ then we denote by $\mathcal{C}^{0}(\mathbb{B})$ or $\mathit{C}(\mathbb{B})$ the vector	 space of complex-valued continuous functions on $\mathbb{B}$ and we denote by $\mathit{C}^{k}(\mathbb{B})$ the set of $k$ times continuously differentiable functions for any $k\in\mathbb{N}$, $k > 0$. The intersection of the spaces $\mathit{C}^{k}(\mathbb{B})$ for all $k\in\mathbb{N}$ is the space $\mathit{C}^{\infty}(\mathbb{B})$ of functions on $\mathbb{B}$ which are differentiable to all orders. It is easy to check that $f \in\mathit{C}^{k}(\mathbb{B})$ if and only if $D^{\alpha\overline{\beta}}f\in\mathit{C}(\mathbb{B})$ for any pair $\left(\alpha,\beta\right)\in\mathbb{N}^{n}\times\mathbb{N}^{n}$ such	that $|\alpha| + |\beta|\leq k$. If $k\in\mathbb{N}$ then the vector space of functions $f$ contained
	in $\mathit{C}^{k}(\mathbb{B})$ whose derivatives $D^{\alpha}f$, $|\alpha|\leq k$, are continuous on $D$ is denoted by $\mathit{C}^{k}(\mathbb{B})$ and we denote by $\mathit{C}^{\infty}(\mathbb{B})$ space of infinitely differentiable functions on $\mathbb{B}$ all of whose derivatives are continuous on $\mathbb{B}$ \cite{Zhu04}.
	
	If $\mathbb{B}\Subset\mathbb{R}^{n}$ and $f\in\mathit{C}^{k}\left(\overline{\mathbb{B}}\right)$, $k\in\mathbb{N}$, then we define the $\mathit{C}^{k}$ norm of $f$ on $\mathbb{B}$ by
	\begin{equation}
	\left\|f\right\|_{k,\mathbb{B}}=\sum_{\alpha\in\mathbb{N}^{n} \\ \left|\alpha\right|\leq k}\underset{x\in\mathbb{B}}{\text{sup}}\left|D^{\alpha}f(x)\right|;
	\end{equation}
	
	Let $\mathbb{B}$ be an open set in $\mathbb{C}^{n}$. A complex-valued function $f$ defined on $\mathbb{B}$ is said to be holomorphic on $\mathbb{B}$ if $f\in\mathit{C}^{1}(\mathbb{B})$ and
	\begin{equation}
	\frac{\partial f}{\partial\overline{z}_{j}}(z)=0
	\end{equation}
	for every $z\in\mathbb{B}$ and $j=1,\cdots, n.$
	The system of partial differential equations (21) is called the \textit{homogeneous	Cauchy–Riemann system}.

	The mapping $F:\mathbb{B}_{n}\to\mathbb{C}^{N}$, where $N$ is a positive integer, is given by $n$ functions as follows:
	\begin{equation*}
	F(z)=(f_{1}(z),\cdots,f_{n}(z)), \ z\in\mathbb{B}.
	\end{equation*}
	We say that $F$ is a holomorphic mapping if each $f_{k}$ is holomorphic in $\mathbb{B}$.
	
	It is clear that any holomorphic mapping $F:\mathbb{B}\to\mathbb{C}^{n}$ has a Taylor type expansion
	\begin{equation*}
	F(z)=\sum a_{\alpha}z^{\alpha},
	\end{equation*}
	where $\alpha = (\alpha_{1} ,\cdots,\alpha_{n})$ is a multi-index of nonnegative integers and each $a_{\alpha}$ belongs to $\mathbb{C}^{n}$.
	
	A mapping $F:\mathbb{B}\to\mathbb{B}$ is said to be bi-holomorphic if
	\begin{itemize}
		\item  $F$ is one-to-one and onto.
		\item  $F$ is holomorphic.
		\item  $F^{-1}$ is holomorphic.
	\end{itemize}

	The automorphism group of $\mathbb{B}$, denoted by $Aut(\mathbb{B})$, consists of all bi-holomorphic mappings of $\mathbb{B}$. It is clear that $Aut(\mathbb{B})$ is a group with composition being the group operation. Conventionally, bi-holomorphic mappings are also called automorphisms.

Let $\mathit{H}^{\infty}(\mathbb{B})$ be the set of all bounded holomorphic functions on $\mathbb{B}.$	The function $||\cdot||_{\infty}:\mathit{H}^{\infty}(\mathbb{B})\to\mathbb{R}_{\geq 0}$ such that
	\begin{equation}\label{key}
	\left\|f\right\|_{\infty}=\mathrm{sup}\left\{\left|f(z)\right|:z\in\mathbb{B}\right\},  \  \forall \ f\in\mathit{H}^{\infty}(\mathbb{B}),
	\end{equation}
	is called the sup norm on $\mathit{H}^{\infty}(\mathbb{B})$, and $\left\|f\right\|_{\infty}$ is called the sup norm of $f$.
	
	The sup norm on $\mathit{H}^{\infty}(\mathbb{B})$ has all the properties of a norm.
	
	A Cauchy sequence in $\mathit{H}^{\infty}$ is a sequence $\left\{f_{n}\right\}_{n\in\mathbb{N}}$, with the following property: for each $\epsilon>0$ there exists an $n(\epsilon)$ such that
	\begin{equation}
	\left\|f_{n}-f_{m}\right\|\leq\epsilon, \ \forall \ n\geq n(\epsilon), \mathrm{and} \ m\geq n(\epsilon).
	\end{equation}
	\begin{thm}
\normalfont
		The algebra of bounded holomorphic functions with the sup norm is a Banach algebra.
	\end{thm}
\subsection{The Space $\mathcal{L}^{p,\infty}$} 	
	Let $\omega$ be a measure space with a (positive) sigma-finite measure $\mu$. The weak	$\mathcal{L}^{p}$ space $\mathcal{L}^{p,\infty}(\mu)$, $0 < p < \infty$, consists of those measurable functions $f$ on $\Omega$ for which
	\begin{equation}
	\left\|f\right\|_{p,\infty}:=\underset{0<\lambda<\infty}{\text{sup}}\lambda\cdot\left(\mu\left(f,\lambda\right)\right)^{1/p}<\infty,
	\end{equation}
	where $\mu\left(f,\lambda\right)=\mu\left\{x:\left|f(x)\right|>\lambda\right\}=\mu\left(\left\{x\in\Omega:\left|f(x)\right|>\lambda\right\}\right).$
	
	Chebyshev's inequality, $$\mu\left(g,\lambda\right)\leq\frac{1}{\lambda}\int_{\Omega}\left|g\right|d\mu,$$ shows that $\mathcal{L}^{p}\subset\mathcal{L}^{p,\infty}$, while the formula
	\begin{equation}
	\int_{\Omega}\left|g\right|^{q}d\mu=\int_{0}^{\infty}\mu\left(g,\lambda\right)d\left(\lambda^{q}\right)
	\end{equation}
	(proved by means of Fubini’s theorem) implies $\mathcal{L}^{p,\infty}\subset\mathcal{L}^{q}$ for $q < p$, if $\mu$ is finite.  The quantity $\left\|\cdot\right\|_{p,\infty}$ is a norm for number $p$, but we have $$\left\|f+g\right\|_{p,\infty}\leq\mathit{C}_{p}\left(\left\|f\right\|_{p,\infty}+\left\|f\right\|_{p,\infty}\right) \ \left(\mathit{C}_{p}=2^{\text{max}(1/p,1)}\right),$$ and hence $\left\|\cdot\right\|_{p,\infty}$ is a (complete) quasinorm. It is interesting, however, that if	$p = 1$, then the space need not be locally convex (if, for example, $\Omega = \left[0, 1\right]$ with Lebesgue measure), although it can be $q-$re-normed for every $q < 1$. For $p > 1$ the
	space is locally convex, and for $p < 1$ it is $p-$convex, i.e., there is an equivalent $p-$norm on it. Hence, the following inequalities holds $$ \mu\left(f_{1}+f_{2},\lambda_{1}+\lambda_{2}\right)\leq\mu\left(f_{1},\lambda_{1}\right)+\left(f_{2},\lambda_{2}\right),$$ and $$\mu\left(f_{1}f_{2},\lambda_{1}\lambda_{2}\right)\leq\mu\left(f_{1},\lambda_{1}\right)+\mu\left(f_{2},\lambda_{2}\right).$$
\subsection{Smooth Manifold}
	Let $U\subset\mathbb{R}^{n}$ and $V\subset\mathbb{R}^{m}$ be open sets, then $f:U\to V$ is called a diffeomorphism or $\mathcal{C}^{\infty}-$diffeomorphism if it is a $\mathcal{C}^{\infty}-$bijection and its inverse $f^{-1}:V\to U$ is $\mathcal{C}^{\infty}$ on $V$.
	
	If $f:U\to V$ is a $\mathcal{C}^{\infty}-$diffeomorphism then for any $a\in U$ $D_{a}f:\mathbb{R}^{n}\to\mathbb{R}^{m}$ is an invertible linear map. In particular $n=m$.
	
	We say that $f:U\to\mathbb{R}^{n}$ is a local diffeomorphism if any $x\in U$ has an open neighbourhood $U_{x}$ such that $$f:U_{x}\to f\left(U_{x}\right)$$ is a diffeomorphism. The inverse function theorem then implies that for any $a\in U$, $D_{a}f$ is invertible.
	
	A pair $\left(M, \mathcal{A}\right)$ where $M$ is a topological space which is separable and second countable, and $\mathcal{A}$ is a collection of continuous maps $\left\{\phi_{\alpha}:U_{\alpha}\to M \ \forall\alpha\in I\right\}$, for open sets $U_{\alpha}\subset\mathbb{R}^{n}$ is a smooth manifold if the following conditions are satisfied:
	\begin{enumerate}
		\item [(1)] $\phi:U_{\alpha}\to\phi\left(U_{\alpha}\right)$ is a homeomorphism, and $$\underset{\alpha}{\bigcup}\phi_{\alpha}\left(U_{\alpha}\right)=M,$$
		\item[(2)] The charts $\left(\phi_{\alpha},U_{\alpha}\right)$ are smoothly compatible i.e. for any $\alpha,\beta\in I$, with $$\phi_{\alpha}\left(U_{\alpha}\right)\cap\phi_{\beta}\left(U_{\beta}\right)\neq\emptyset,$$ $$\phi_{\alpha}^{-1}\circ\phi_{\beta}:\phi_{\beta}^{-1}\left(\phi_{\alpha}\left(U_{\alpha}\right)\cap\phi_{\beta}\left(U_{\beta}\right)\right)\to\phi_{\alpha}^{-1}\left(\phi_{\alpha}\left(U_{\alpha}\right)\cap\phi_{\beta}\left(U_{\beta}\right)\right)$$
		is a diffeomorphism. The above condition makes $\mathcal{A}$ a smooth atlas.
		\item[(3)] $\mathcal{A}$ is a maximal smooth atlas i.e. there is strictly no larger smooth atlas containing $\mathcal{A}$.
		\end{enumerate}	
\begin{definition}
\normalfont
	A Lie group is a group $G$ which is a differentiable manifold and such that multiplication and inversion are smooth  maps. By multiplication we mean that $\left(g,h\right)\mapsto gh:G\times G\to G$ is smooth and by inversion we mean that $g\mapsto g^{-1}:G\to G$ are smooth.
	
	Let $G$ and $H$ be two Lie groups. Then $f:G\to H$ is a Lie group morphism if it is smooth and a group morphism. If $H\subset G$, then $H$ is a Lie subgroup of $G$ if it is at the same time a subgroup and submanifold of $G$.
	
	A Lie algebra on a field $\mathbb{K}$ is a $\mathbb{K}-$vector space $V$ endowed with an antisymmetric $\mathbb{K}-$bilinear map $\left[\cdot,\cdot\right]:V\times V\to V$ which satisfy the Jacobi identity $$\left[\left[x,y\right],z\right]+\left[\left[y,z\right],x\right]+\left[\left[z,x\right],y\right]=0.$$ By antisymmetric we mean that $$\left[x,y\right]=-\left[y,x\right]$$ and $\left[\cdot,\cdot\right]$ is called the Lie bracket.
		
	A Lie subalgebra of $V$ is a linear subspace of $V$ which is a linear subspace of $V$ which is  stable under $\left[\cdot,\cdot\right].$
		
	A lie algebra morphism is a linear map $A$ from $\left(V,\left[\cdot,\cdot\right]_{V}\right)$ to $\left(V^{'},\left[\cdot,\cdot\right]_{V^{'}}\right)$ such that $$\left[A(x),A(y)\right]_{V^{'}}=A\left(\left[x,y\right]_{V}\right).$$ A Lie algebra isomorphism is a bijective Lie algebra morphism, it is then automatic that its inverse is a Lie algebra morphism.
\end{definition}
	\section{Algebra of Bergman Spaces}
Here, we shall study the properties of the Bergman spaces and the Toeplitz operator algebra defined on them.
\begin{definition}
\normalfont
	For $\alpha>0$ the weighted Lebesgue measure $d\mathcal{A}_{\alpha}$ is defined by $$d\mathcal{A}_{\alpha}(z)=C_{\alpha}\left(1-|z|^{2}\right)^{\alpha}d\mathcal{A}(z)$$ where $$c_{\alpha}=\frac{\Gamma\left(n+\alpha+1\right)}{n!\Gamma\left(\alpha+1\right)},$$  is a normalizer so that $d\mathcal{A}_{\alpha}$ forms a probability measure on $\mathbb{B}_{n}$. For $\alpha>-1$ and $p>0$ the weighted Bergman space $\mathcal{A}_{\alpha}^{p}$ consists of holomorphic functions $f$ in $\mathit{L}^{p}(\mathbb{B},d\mathcal{A})$, that is,
	\begin{equation*}
	\mathcal{A}_{p}=\mathit{L}^{p}(\mathbb{B},d\mathcal{A})\cap\mathit{H}(\mathbb{B}).
	\end{equation*}
\end{definition}
	It is clear that $\mathcal{A}_{\alpha}^{p}$ is a linear subspace of $\mathit{L}^{p}(\mathbb{B},d\mathcal{A})$. When the weight $\alpha =0$, we simply write $\mathcal{A}^{p}$ for $\mathcal{A}_{\alpha}^{p}$. These give the standard (unweighted) Bergman spaces.

The Bergman norm is a quasinorm  given by
	\begin{equation}
	\left\|f\right\|=\left\|f\right\|_{\mathcal{A}^{p}}=\left(2\int_{0}^{1}\mathit{I}_{p}(r,f)rdr\right)^{\frac{1}{p}},
	\end{equation}
	where \begin{equation}
	\mathit{I}_{p}(r,f)=\left(\frac{1}{2\pi}\int_{-\pi}^{\pi}\left|f(re^{i\theta})\right|^{p}d\theta\right)^{\frac{1}{p}}.
	\end{equation}

	It is clear that $\mathcal{A}^{p}$ is also a linear subspace $\mathit{L}^{p}(\mathbb{B},d\mathcal{A})$.
	\begin{prop}\cite{Pavlovic}
\normalfont
	Let $p\in(0,\infty)$. Then, the following hold:
	\begin{enumerate}
		\item For every $z\in\mathbb{B}$ the functional $f\mapsto f(z)$ is continuous on $\mathcal{A}^{p}$;
		\item The space $\mathcal{A}^{p}$ is complete;
		\item If $f\in\mathcal{A}^{p}$, and $f_{n}(z)=f(nz)$, then $f_{n}\in\mathcal{A}^{p}$ and $\left\|f-f_{n}\right\|\to0(n\to\infty)$;
		\item The set of all polynomials is a dense subset $\mathcal{A}^{p}$.
	\end{enumerate}
	\end{prop}
\begin{prop}
\normalfont
		The space $\mathit{L}^{p}(\mathbb{B}, d\mathcal{A})$ is a Banach algebra with convolution operation defined as
	\begin{equation}
	\left\|f* g\right\|^{p}=2\int_{0}^{1}\mathit{I}_{p}\left(r,f* g\right)rdr
	\end{equation}
	for all $f,g\in\mathit{L}^{p}(\mathbb{B}, d\mathcal{A})$
	\end{prop}
	\paragraph*{Proof.} Now,
	\begin{equation*}
	\left\|f* g\right\|^{p}=2\int_{0}^{1}\mathit{I}_{p}\left(r,f* g\right)rdr.
	\end{equation*}
	To this end,
	\begin{flalign*}
	\mathit{I}_{p}\left(r,f*g\right)&=\frac{1}{2\pi}\int_{-\pi}^{\pi}\left|f* g\right|^{p}\left(re^{i\theta}\right)d\theta\\
	&=\frac{1}{2\pi}\int_{-\pi}^{\pi}\left|\frac{1}{2\pi}\int_{-\pi}^{\pi}f(x)g(x-y)dy\right|^{p}\left(re^{i\theta}\right)d\theta\\
	&\leq\frac{1}{2\pi}\int_{-\pi}^{\pi}\frac{1}{2\pi}\int_{-\pi}^{\pi}\left|f(x)g(x-y)dy\right|^{p}\left(re^{i\theta}\right)d\theta.
	\end{flalign*}
	Letting $\theta=x-y$ so that $x=\theta+y$ and $d\theta=dx,$ we obtain
	\begin{equation*}
	re^{i\theta}=re^{i(x-y)}=re^{ix-iy}=re^{ix}\cdot re^{-iy}.
	\end{equation*}
	Hence,
	\begin{flalign*}
	\left\|f* g\right\|^{p}&=\frac{1}{2\pi}\int_{-\pi}^{\pi}\frac{1}{2\pi}\int_{-\pi}^{\pi}\left|f(\theta +y)g(\theta)\right|^{p}dy\left(re^{ix}\right)\left(re^{-iy}\right)dx\\
	&\leq\frac{1}{2\pi}\int_{-\pi}^{\pi}\frac{1}{2\pi}\int_{-\pi}^{\pi}\left|f(\theta +y)\right|^{p}\left(re^{-iy}\right)dy\left|g(z)\right|^{p}\left(re^{ix}\right)dx\\
	&=\frac{1}{2\pi}\int_{-\pi}^{\pi}\left|f(\theta +y)\right|^{p}re^{-iy}dy\frac{1}{2\pi}\int_{-\pi}^{\pi}\left|g(\theta)\right|^{p}\left(re^{ix}\right)dx\\
	&=\frac{1}{2\pi}\int_{-\pi}^{\pi}\left|f(\theta -y)\right|^{p}re^{iy}dy\frac{1}{2\pi}\int_{-\pi}^{\pi}\left|g(\theta)\right|^{p}\left(re^{ix}\right)dx\\
	&=\mathit{I}_{p}\left(r,f\right)\cdot\mathit{I}_{p}\left(r,g\right)
	\end{flalign*}
	For this reason,
	\begin{equation*}
	\left\|f* g\right\|_{\mathcal{A}^{p}}\leq\left\|f\right\|_{\mathcal{A}^{p}}\left\|g\right\|_{\mathcal{A}^{p}}. \ \ \ \ \square
	\end{equation*}
\begin{rem}
\normalfont
	In the special case when $p=2, \ \mathit{L}^{2}(\mathbb{B},d\mathcal{A})$ is a Hilbert space whose inner product is denoted by
	\begin{equation}
	\left\langle f,g\right\rangle_{\mathcal{A}^{p}}=\int_{\mathbb{B}}f(z)\overline{g(z)}d\mathcal{A}(z)
	\end{equation}
	and
	\begin{flalign*}
	\left|\left\langle f,g\right\rangle\right|_{\mathcal{A}^{p}}&\leq\left(\int_{\mathbb{B}}\left|f(z)\right|^{2}d\mathcal{A}(z)\right)^{1/2}\left(\int_{\mathbb{B}}\left|g(z)\right|^{2}d\mathcal{A}(z)\right)^{1/2}\\
	&=\left\|f\right\|_{\mathcal{A}^{p}}\left\|g\right\|_{\mathcal{A}^{p}}
	\end{flalign*}
\end{rem}
	\begin{thm}
\normalfont
	Suppose $x$ and $y$ are complex numbers, $\varphi$ and  $\phi$ are bounded functions on $\mathbb{B}$; then
	\begin{enumerate}
	\item $T_{x\varphi}+T_{y\phi}=xT_{\varphi}+yT_{\phi}$;
	\item $T_{\overline{\varphi}}=T_{\varphi}^{*}$;
	\item $T_{\varphi}\geq0 \ \mathrm{if} \ \varphi\geq 0$;\\
	Moreover, if $\varphi\in\mathit{H}^{\infty}$, then
	\item $T_{\phi}T_{\varphi}=T_{\phi\varphi}$;
	\item $T_{\overline{\varphi}}T_{\phi}=T_{\overline{\varphi}\phi}$.
	\end{enumerate}
	\end{thm}
	\paragraph*{Proof:}
	Given that, $\varphi,\phi\in\mathit{L}^{\infty}(\mathbb{B})$, then for any $x,y\in\mathbb{C}$, (1) implies
	\begin{flalign*}
	T_{x\varphi}(f)+T_{y\varphi}(f)&=P(x\varphi f)+P(y\phi f)\\
	&=\int_{\mathbb{B}}K(z,u)(x\varphi)(u)f(u)d\mathcal{A}(u)+\int_{\mathbb{B}}K(z,u)(y\varphi)(u)f(u)d\mathcal{A}(u)\\
	 &=\int_{\mathbb{B}}\frac{(x\varphi)(u)f(u)}{(1-z\overline{u})^{2}}d\mathcal{A}(u)+\int_{\mathbb{B}}\frac{(y\varphi)(u)f(u)}{(1-z\overline{u})^{2}}d\mathcal{A}(u)\\
	 &=\int_{\mathbb{B}}\frac{x\varphi(u)f(u)}{(1-z\overline{u})^{2}}d\mathcal{A}(u)+\int_{\mathbb{B}}\frac{y\varphi(u)f(u)}{(1-z\overline{u})^{2}}d\mathcal{A}(u)\\
	&=	 x\int_{\mathbb{B}}\frac{\varphi(u)f(u)}{(1-z\overline{u})^{2}}d\mathcal{A}(u)+y\int_{\mathbb{B}}\frac{\varphi(u)f(u)}{(1-z\overline{u})^{2}}d\mathcal{A}(u)\\
	&=x\int_{\mathbb{B}}K(z,u)\varphi(u)f(u)d\mathcal{A}(u)+y\int_{\mathbb{B}}K(z,u)\varphi(u)f(u)d\mathcal{A}(u)\\
	&=xP(\varphi f)+yP(\phi f)\\
	&=xT_{\varphi}(f)+yT_{\phi}(f).
	\end{flalign*}
	(2) implies;
	\begin{flalign*}
	T_{\overline{\varphi}}(f)&=P(\overline{\varphi}f)\\
	&=\int_{\mathbb{B}}K(z,u)\overline{\varphi}(u)f(u)d\mathcal{A}(u)\\
	&=\int_{\mathbb{B}}\frac{\overline{\varphi}(u)f(u)}{(1-z\overline{u})^{2}}d\mathcal{A}(u)\\
	&=\int_{\mathbb{B}}\frac{\left(\varphi(u)f(u)\right)^{*}}{(1-z\overline{u})^{2}}d\mathcal{A}(u)\\
	&=\int_{\mathbb{B}}K(z,u)\left(\varphi(u)f(u)\right)^{*}d\mathcal{A}(u)\\
	&=P(\varphi f)^{*}\\
	&=T_{\varphi}^{*}(f).
	\end{flalign*}
	(3.) Since there are non zero divisors in the set of  all Toeplitz operators for all $\varphi\in\mathit{L}^{\infty}(\mathbb{B})$, then
	\begin{equation*}
	T_{\varphi}(f)=0
	\end{equation*}
	implies
	\begin{flalign*}
	T_{\varphi}(f)&=P(\varphi f)(z)\\
	&=\int_{\mathbb{B}}K(z,u)\varphi(u)f(u)d\mathcal{A}(u)\\
	&=\int_{\mathbb{B}}\frac{\varphi(u)f(u)}{(1-z\overline{u})^{2}}d\mathcal{A}(u)\\
	&=0.
	\end{flalign*}
	Since $f(u)\neq 0$, $\varphi(u)$ must be equal to zero. Also, $T_{\varphi}(f)>0$ implies $$T_{\varphi}(f)=P(\varphi f)>0 \ \ \implies \ \ \varphi> 0.$$ Hence, $$T_{\varphi}(f)\geq 0 \ \implies \ P(\varphi f)\geq 0 \ \implies \ \varphi\geq 0.$$

	(4) If $\varphi, \phi\in\mathit{H}^{\infty}$, then $\varphi\mathcal{A}^{2}\subseteq\mathcal{A}^{2}$ and,
	\begin{flalign*}
	T_{\phi}T_{\varphi}(f)&=P(\phi f)P(\varphi f)\\
	&=\left\langle\int_{\mathbb{B}}K(z,u)\phi(u)f(u)d\mathcal{A}(u),\int_{\mathbb{B}}K(z,u)\varphi(u)f(u)d\mathcal{A}(u)\right\rangle\\
	&=\int_{\mathbb{B}}K(z,u)\left\langle\phi,\varphi\right\rangle f(u)d\mathcal{A}(u)\\
	&=\int_{\mathbb{B}}\frac{\left\langle\phi,\varphi\right\rangle f(u)}{(1-z\overline{u})^{2}}d\mathcal{A}(u)\\
	&=\int_{\mathbb{B}}\frac{\phi\varphi f(u)}{(1-z\overline{u})^{2}}d\mathcal{A}(u)\\
	&=\int_{\mathbb{B}}K(z,u)\phi(u)\varphi(u)f(u)d\mathcal{A}(u)\\
	&=P(\phi\varphi f)\\
	&=T_{\phi\varphi}(f).
	\end{flalign*}
	(5) follows from (2) and (4) that is, by taking the adjoint,
	\begin{flalign*}
	T_{\overline{\varphi}}(f)T_{\phi}(f)&=T_{\varphi}^{*}(f)T_{\overline{\phi}}^{*}(f)\\
	&=\left(T_{\varphi}(f)T_{\overline{\phi}}(f)\right)^{*}\\
	&=\left(T_{\varphi\overline{\phi}}(f)\right)^{*}\\
	&=T_{\overline{\varphi}\phi}(f)
	\end{flalign*}
	and therefore, the proof is complete. \ \ \ \ $\square$
\section{Algebra of Bloch Spaces}
	\begin{definition}
\normalfont
	For $0\leq\alpha<\infty$, let $\mathcal{H}_{\alpha}^{\infty}$ be the space of holomorphic functions $f\in\mathcal{H}(\mathbb{B})$ satisfying
	\begin{equation}
	\underset{z\in\mathbb{B}}{\mathrm{sup}}\left(1-\Vert z\Vert^{2}\right)^{\alpha}\left|f(z)\right|<\infty.
	\end{equation}
	We abbreviate $\mathcal{H}^{\infty}=\mathcal{H}_{1}^{\infty}$ for $\alpha=1$.
	\end{definition}
	
	The classical $\alpha-$Bloch space, denoted as $\mathcal{B}^{\alpha}$ is the space of holomorphic functions $F : \mathbb{B}\to\mathbb{C}$ satisfying
	\begin{equation}
	 \left\|f\right\|_{\mathcal{B}^{\alpha}(\mathbb{B})}\:=\underset{z\in\mathbb{B}}{\mathrm{sup}}\left(1-\left\|z\right\|^{2}\right)^{\alpha}\left|f'(z)\right|<+\infty.
	\end{equation}
	Now we introduce four semi-norms of the Bloch type space (see \cite{Xu}) for $f\in\mathcal{H}(\mathbb{B})$ in what follows. To do this, let
	\begin{equation}
	 \left\|f\right\|_{\mathcal{B},\alpha}:=\underset{z\in\mathbb{B}}{\mathrm{sup}}\left(1-\left\|z\right\|^{2}\right)^{\alpha}\left|\partial^{m}f(z)\right|,
	\end{equation}
	\begin{equation}
	\left\|f\right\|_{\mathcal{R},\alpha}:=\underset{z\in\mathbb{B}}{\mathrm{sup}}\left(1-\Vert z\Vert^{2}\right)^{\alpha}\left|\mathcal{R}f(z)\right|,
	\end{equation}
	\begin{equation}
	\left\|f\right\|_{weak,\alpha}:=\underset{z\in\mathbb{T}}{\mathrm{sup}}\left\|f_{y}\right\|_{\mathcal{B}^{\alpha}(\mathbb{B})},
	\end{equation}
	where $\mathcal{R}f(z)=\left\langle\partial^{m}f(z),z\right\rangle, \ f_{y}(z)=f(zy)$ for $z\in\mathbb{B}, \  |z|<1$, for each $y\in\mathbb{B}$ with norm $||y||=1$. We note that $zf_{y}^{\prime}(z)=\mathcal{R}f(zy)$. The Mobius transforms of $\mathbb{B}$ are holomorphic mappings $\varphi_{a}$, $a\in\mathbb{B}$, given by
	\begin{equation}
	\varphi_{a}(z)=\left(\mathit{P}_{a}+s_{a}\mathit{Q}_{a}\right)\left(m_{a}(z)\right),
	\end{equation}
	where $s_{a}=\sqrt{1-\Vert a\Vert^{2}}, \ \mathit{P}_{a}(z)=\frac{\left\langle z,a\right\rangle}{\left\langle a,a\right\rangle}a, \ \mathit{Q}_{a}=\mathit{I}-\mathit{P}_{a}$ and $m_{a}(z)=\frac{a-z}{1-\left\langle z,a\right\rangle}$.
	Define
	\begin{equation}
	 \left\|f\right\|_{\tilde{\mathcal{B}},\alpha}:=\underset{z\in\mathbb{B}}{\mathrm{sup}}\left(1-\left\|z\right\|^{2}\right)^{\alpha-1}\left|\tilde{\mathbf{\nabla}}f(z)\right|,
	\end{equation}
	where $\tilde{\mathbf{\nabla}}=\partial^{m}f\circ\varphi_{z}(0)$ with $\varphi_{z}\in\mathrm{Aut}(\mathbb{B})$.
	We note that, by Lemma 3.5 of \cite{Blasco}
	\begin{equation}
	 \left\|\tilde{\mathbf{\nabla}}f(z)\right\|=\underset{w\neq0}{\mathrm{sup}}\frac{(1-||z||^{2})\left|\partial^{m}f(z)(w)\right|}{\sqrt{(1-||z||^{2})||w||^{2}+\left|\left\langle w,z\right\rangle\right|^{2}}}.
	\end{equation}
	Hence, we have
	\begin{equation}
	 \left\|f\right\|_{\tilde{\mathcal{B}},\alpha}=\underset{x\in\mathbb{B}}{\mathrm{sup}}\underset{w\neq0}{\mathrm{sup}}\frac{(1-||z||^{2})^{\alpha}\left|\partial^{m}f(z)(w)\right|}{\sqrt{(1-||z||^{2})||w||^{2}+\left|\left\langle w,z\right\rangle\right|^{2}}}.
	\end{equation}
	Hence for $\alpha>0$, we have
	\begin{equation}
	\mathcal{B}^{\alpha}=\left\{f\in\mathcal{H}(\mathbb{B}):\left\|f\right\|_{\tilde{\mathcal{B}},\alpha}<+\infty\right\}.
	\end{equation}
	\begin{prop}
\normalfont
	Equipped with the norm $\left\|f\right\|_{\alpha}=|f(0)|+\left\|f\right\|_{\mathcal{B},\alpha}$ for $f\in\mathcal{B}^{\alpha}$, the Bloch type space $\mathcal{B}^{\alpha}$ becomes	a Banach space and hence, a Banach algebra.
	\end{prop}
	\paragraph*{Proof.} For any $f,g\in\mathcal{B}^{\alpha}$ we have,
	\begin{flalign*}
		\Vert fg\Vert_{\alpha}&=\Vert fg(0)\Vert+\underset{z\in\mathbb{B}}{\mathrm{sup}}(1-||z||^{2})^{\alpha}\left\|(fg)^{\prime}(z)\right\|_{\alpha}\\
		&=\Vert fg(0)\Vert+\underset{z\in\mathbb{B}}{\mathrm{sup}}(1-||z||^{2})^{\alpha}\left\|(f^{\prime}g)(z)+(fg^{\prime}(z))\right\|_{\alpha}\\
		 &\leq\left\|f(0)\right\|\left\|g(0)\right\|+\underset{z\in\mathbb{B}}{\mathrm{sup}}(1-||z||^{2})^{\alpha}\left(\left\|f^{\prime}g(z)\right\|_{\alpha}+\left\|fg^{\prime}(z)\right\|_{\alpha}\right)\\
		 &=\left\|f(0)\right\|\left\|g(0)\right\|+\underset{z\in\mathbb{B}}{\mathrm{sup}}(1-||z||^{2})^{\alpha}\left\|f^{\prime}g(z)\right\|_{\alpha}+\underset{z\in\mathbb{B}}{\mathrm{sup}}(1-||z||^{2})^{\alpha}\left\|fg^{\prime}(z)\right\|_{\alpha}\\
		 &\leq\left\|f(0)\right\|\left\|g(0)\right\|+\underset{z\in\mathbb{B}}{\mathrm{sup}}(1-||z||^{2})^{\alpha}\left\|f^{\prime}(z)\right\|_{\alpha}\left\|g(z)\right\|_{\alpha}\\
		&+\underset{z\in\mathbb{B}}{\mathrm{sup}}(1-||z||^{2})^{\alpha}\left\|f(z)\right\|_{\alpha}\left\|g^{\prime}(z)\right\|_{\alpha}\\
		 &=|f(0)||g(0)|+\left\|f^{\prime}(z)\right\|_{\alpha}\left\|g(z)\right\|_{\alpha}+\left\|f(z)\right\|_{\alpha}\left\|g^{\prime}(z)\right\|_{\alpha}.\;\;\;\;\;\;\;\;\;\;\;\;\;\; \Box
	\end{flalign*}
	\begin{definition}[The Little Bloch Space]
\normalfont
	We denote the class of Bloch functions defined on $\mathbb{B}$ by $\mathcal{B}(\mathbb{B})$. The Bloch space is not separable but there exists a separable subspace of the Bloch space known as the little Bloch space.
	
	We let $\mathcal{C}(\overline{\mathbb{B}_{n}})$ be the space of continuous function on the closed unit ball, and $\mathcal{C}_{0}(\mathbb{B}_{n})$ be closed subspace of $\mathcal{C}(\overline{\mathbb{B}_{n}})$ consisting of those functions that vanish on the boundary $\mathbb{T}_{n}$.
	
	The little Bloch space $\mathcal{B}_{0}(\mathbb{B})$ is the subspace of $\mathcal{B}(\mathbb{B})$ given by those functions for which
	\begin{equation*}
	\underset{|z|\to 1-}{\text{lim}}\left|\tilde{\mathbf{\nabla}}f(z)\right|=0.
	\end{equation*}
	Since $\left|\tilde{\mathbf{\nabla}}f(z)\right|$ is dense in $\mathbb{B}_{n}$, the above condition means that the function $\left|\tilde{\mathbf{\nabla}}f(z)\right|$ belongs to $\mathcal{C}_{0}(\mathbb{B}_{n})$.
	\end{definition}
	
	Now, if $f\in\mathcal{B}(\mathbb{B})$ then  $x^{*}f\in\mathcal{B}$ for all $x^{*}\in\mathbb{B}$. And, interchanging the suprema, we have that
	\begin{equation*}
	\left\|f\right\|_{\mathcal{B}(\mathbb{B})}\thickapprox\underset{\Vert x^{*}\Vert=1}{\mathrm{sup}}\Vert x^{*}f\Vert_{\mathcal{B}}
	\end{equation*}
	where $x^{*}f(z)=\left\langle f(z),x\right\rangle$.
	
	The Bloch space possesses the following properties.
	\begin{definition}[Pointwise Multiplier on Bloch Spaces]
\normalfont
	From definition 1.8 number (5), a function $f$ is called a pointwise multiplier of a space $\mathcal{B}\left(\mathbb{B}\right)$ if for every $g \in \mathcal{B}\left(\mathbb{B}\right)$ the	pointwise product $fg$ also belongs to $\mathcal{B}\left(\mathbb{B}\right)$. Thus, we denote a pointwise multiplier $f$ of a space $\mathcal{B}\left(\mathbb{B}\right)$ by $f\mathcal{B}\left(\mathbb{B}\right)\subset \mathcal{B}\left(\mathbb{B}\right)$.
	
	Here, we consider the pointwise multiplier algebra of the Bloch space. Throughout this section we use the following norm on $\mathcal{B}$: $$\left\|g\right\|=\left|g(0)\right|+\text{sup}\left\{\left(1-|z|^{2}\right)\left|\nabla g(z)\right|:z\in\mathbb{B}_{n}\right\}, \ g\in\mathcal{B}.$$
	\end{definition}
	\begin{prop}
\normalfont
	For all $f,g\in\mathcal{B}(\mathbb{B})$ and $\varphi\in\mathit{H}^{\infty}(\mathbb{B})$, the following properties hold.
	\begin{enumerate}
		\item [(i).] $\left\|\mathcal{M}_{a\varphi}(f)\right\|=|a|\left\|\mathcal{M}_{\varphi}(f)\right\|$,
		\item[(ii)] $\left\|\mathcal{M}_{a\varphi}(f)\right\|+\left\|\mathcal{M}_{b\varphi}(f)\right\|=|a|\left\|\mathcal{M}_{\varphi}(f)\right\|+|b|\left\|\mathcal{M}_{\varphi}(f)\right\|$ and \item[(iii).] $\left\|\mathcal{M}_{\varphi_{1}\varphi_{2}}(f)\right\|\leq\left\|\varphi_{1}\right\|_{\infty}\left\|\varphi_{2}\right\|_{\infty}$,
	\end{enumerate}
	for all $\varphi_{1},\varphi_{2}\in\mathit{H}^{\infty}$ and $a,b$ are constants. In particular, $\mathcal{M}_{\varphi}$ is a Banach algebra.
	\end{prop}
	\paragraph*{Proof.} Given that $$\left\|f(z)\right\|=\left|f(0)\right|+\left\|f(z)\right\|$$ where $\left\|f(z)\right\|=\underset{z\in\mathbb{B}}{\text{sup}}\left(1+|z|^{2}\right)\left|\nabla f(z)\right|$ for $f\in\mathcal{B}$, then (i) implies
	\begin{flalign*}
	\left\|\mathcal{M}_{a\varphi}(f)\right\|&=\left|(a\varphi)f(0)\right|+\left\|(a\varphi)f(z)\right\|\\
	&\leq\left|a\right|\left\|\varphi f(0)\right\|+\left|a\right|\left\|\varphi f(z)\right\|\\
	&\leq\left|a\right|\left\|\varphi\right\|_{\infty}\left(\left|f(0)\right|+\left\|f(z)\right\|\right)\\
	&=|a|\left\|\varphi\right\|_{\infty}\left\|f(z)\right\|_{\mathcal{B}}\\
	&|a|\left\|\mathcal{M}_{\varphi}(f)\right\|.
	\end{flalign*}
	(ii) follows immediately from (i).
	And now (iii) can be shown  by applying the close graph theorem as follows;
	\begin{flalign*}
	 \left\|\mathcal{M}_{\varphi_{1}\varphi_{2}}(f)\right\|&=\left|\left(\varphi_{1}\varphi_{2}f(0)\right)\right|+\left\|\varphi_{1}\varphi_{2}f(z)\right\|\\
	&\leq\left\|\varphi_{1}\varphi_{2}\right\|\left(\left|f(0)\right|+\left\|f(z)\right\|\right)\\
	&\leq\left\|\varphi_{1}\varphi_{2}\right\|_{\infty}\left\|f(z)\right\|_{\mathcal{B}}\\
	 &\leq\underset{f\in\mathcal{B}}{\text{sup}}\frac{\left\|\varphi_{1}\varphi_{2}\right\|\left\|f(z)\right\|_{\mathcal{B}}}{\left\|f(z)\right\|_{\mathcal{B}}}\\
	&=\left\|\varphi_{1}\varphi_{2}\right\|_{\infty}\\
	&\leq\left\|\varphi_{1}\right\|_{\infty}\left\|\varphi_{2}\right\|_{\infty}. \ \ \ \ \square
	\end{flalign*}	
	\section{Algebra of Hardy Spaces}
	\begin{definition}
\normalfont
	For $0<p<\infty$ and $\mathit{T}_{n}$ a sphere of radius $n$ which is the boundary of $\mathbb{B}_{n}$, the Hardy space $\mathit{H}^{p}$ consists of holomorphic functions $f$ in $\mathbb{B}_{n}$ such that
	\begin{equation}
	\left\|f\right\|_{p}^{p}=\underset{0<r<1}{\mathrm{sup}}\int_{\mathbb{T}_{n}}\left|f(re^{i\theta})\right|^{p}dm(\theta)<\infty.
	\end{equation}
	\end{definition}
	
	Just as the Bergman kernel plays an essential role in the study of Bergman spaces, two integral kernels are fundamental in the theory of Hardy spaces; they are the Cauchy-Szego kernel,
	\begin{equation}\label{key}
	\mathbf{C}_{s}(z,\xi)=\frac{1}{(1-\left\langle z,\xi\right\rangle^{2})^{n}},
	\end{equation}
	and the (invariant) Poisson kernel,
	\begin{equation}\label{key}
	\mathit{P}(z,\xi)=\frac{(1-|z|^{2})^{n}}{\left|1-\left\langle z,\xi\right\rangle\right|^{2n}}.
	\end{equation}
	We note that, the Poisson kernel here is different from the associated Poisson kernel
	when $\mathbb{B}_{n}$ is thought of as the unit ball in $\mathbb{R}^{2n}$, unless $n=1$.
	Hence, if $f$ belongs to the ball algebra, then
	\begin{gather}
	f(z)=\int_{\mathbb{T}_{n}}\mathbf{C}_{s}(z,\xi)f(\xi)d\mu(\xi) \\
	f(z)=\int_{\mathbb{T}_{n}}\mathit{P}(z,\xi)f(\xi)d\mu(\xi) \\
	\left|f(z)\right|^{p}\leq\int_{\mathbb{T}_{n}}\mathit{P}(z,\xi)\left|f(\xi)\right|^{p}d\mu(\xi)
	\end{gather}		
	for all $z\in\mathbb{B}_{n}$ and $0<p<\infty$.
	\begin{thm}
\normalfont
		If $f$ and $g$ are bounded on $\mathbb{B}$ and $a, b$ are complex numbers, then
	\begin{enumerate}
		\item $T_{af+bg}=aT_{f}+bT_{g}$;
		\item $T_{\overline{f}}=T_{f}^{*}$,  where $T^{*}$ is the adjoint of $T$;
		
		If further $f$ is in $\mathit{H}^{\infty}$, then
		\item $T_{g}T_{f}=T_{gf}$;
		\item $T_{\overline{f}}T_{g}=T_{\overline{f}g}$
	\end{enumerate}
	\end{thm}
	\paragraph{Proof:}
	These properties follow easily from the definition of Toeplitz operators and simple calculations with the inner product in $\mathit{H}^{2}$ shows that:
	for all $f,g\in\mathit{H}^{2}$ and $a,b\in\mathbb{B}$, (1) implies
	\begin{equation*}
	\begin{split}
	T_{af}+T_{bg}&=\mathbf{C}_{s}(af)+\mathbf{C}_{s}(bg)\\
	&=\int_{\mathbb{B}}\frac{af(\xi)}{(1-\left\langle z,\xi\right\rangle)^{n}}d\beta(\xi)+\int_{\mathbb{B}}\frac{bg(\xi)}{(1-\left\langle z,\xi\right\rangle)^{n}}d\beta(\xi)\\
	&=a\int_{\mathbb{B}}\frac{f(\xi)}{(1-\left\langle z,\xi\right\rangle)^{n}}d\beta(\xi)+b\int_{\mathbb{B}}\frac{g(\xi)}{(1-\left\langle z,\xi\right\rangle)^{n}}d\beta(\xi)\\
	&=a\mathbf{C}_{s}(f)+b\mathbf{C}_{s}(g)\\
	&=aT_{f}+bT_{g}.
	\end{split}
	\end{equation*}
	(2) $\implies$
	\begin{equation*}
	\begin{split}
	T_{\overline{f}}&=\mathbf{C}_{s}(\overline{f})\\
	&=\int_{\mathbb{B}}\frac{\overline{f}(\xi)}{(1-\left\langle z,\xi\right\rangle)^{n}}d\beta(\xi)\\
	&=\int_{\mathbb{B}}\frac{f^{*}(\xi)}{(1-\left\langle z,\xi\right\rangle)^{n}}d\beta(\xi)\\
	&=\mathbf{C}_{s}(f)^{*}\\
	&=T_{f}^{*}.
	\end{split}
	\end{equation*}
	Now, given that $f,g\in\mathit{H}^{\infty}$, then (3) implies
	\begin{equation*}
	\begin{split}
	T_{g}T_{f}&=\mathbf{C}_{s}(g)\mathbf{C}_{s}(f)\\
	&=\int_{\mathbb{B}}\frac{g(\xi)}{(1-\left\langle z,\xi\right\rangle)^{n}}d\beta(\xi)\int_{\mathbb{B}}\frac{f(\xi)}{(1-\left\langle z,\xi\right\rangle)^{n}}d\beta(\xi)\\
	&=\int_{\mathbb{B}}\frac{1}{{(1-\left\langle z,\xi\right\rangle)^{n}}}d\beta(\xi)\int_{\mathbb{B}}\frac{g(\xi)f(\xi)}{(1-\left\langle z,\xi\right\rangle)^{n}}d\beta(\xi)\\&=\mathbf{C}_{1}\int_{\mathbb{B}}\frac{f(\xi)g(\xi)}{(1-\left\langle z,\xi\right\rangle)^{n}}d\beta(\xi)\\
	&=\mathbf{C}_{s_1}\mathbf{C}_{s_2}(gf)\\
	&=\mathbf{C}_{s}(gf)\ \ \left(\mathrm{taken} \ \mathbf{C}_{s_1}\mathbf{C}_{s_2}=\mathbf{C}_{s}\right)\\
	&=T_{gf}.
	\end{split}
	\end{equation*}
	Combining (2) and (3), (4) implies
	\begin{equation*}
	\begin{split}
	T_{\overline{f}}T_{g}&=\mathbf{C}_{s}(\overline{f})\mathbf{C}_{s}(g)\\
	&=\mathbf{C}_{s}(\overline{f}g)\\
	&=T_{\overline{f}g}.\ \ \ \ \square
	\end{split}
	\end{equation*}
	\section{Algebra of Spherical Functions.}
	In this section, we consider the algebra of spherical functions. They are the Gelfand pair analogue of characters on locally compact Abelian groups. Let $G$ be a locally compact group and $K$ be a compact subgroup. In the foremost, we do not require $\left(G,K\right)$ to be a Gelfand pair.
	\begin{definition}\cite{Wolf},\cite{Shinya}
\normalfont	
	A \textbf{spherical measure} for $\left(G,K\right)$ is a non zero Radon measure $\mu$ on $G$ such that
	\begin{enumerate}
		\item[(1).] $\mu$ is $K$ bi-invariant i.e. $\mu(k_{1}Ek_{2}^{-1})=\mu(E)$ for every Borel set $E\subset G$, and
		\item[(2).] $f\mapsto\mu(f)=\displaystyle\int_{G}f(g)d\mu(g)$ is an algebra homomorphism $C_{c}(K\setminus G/K)\to\mathbb{C}$, i.e.
		\begin{equation}\label{key}
		\mu(f_{1}* f_{2})=\mu(f_{1})\mu(f_{2}).
		\end{equation}
		In other words, $\mu$ is multiplicative linear functional on $C_{c}(K\setminus G/K)$.
	\end{enumerate}
	\end{definition}
	\begin{definition}\cite{Helgason}
\normalfont
	A \textbf{Spherical function} for $(G,K)$ is a continuous function $\phi:G\to\mathbb{C}$ such that the integral $\phi(f)$ defined by
	\begin{equation}\label{key}
	\phi(f)=\int_{G}f(x)\phi(x^{-1})d\mu_{G}(x)
	\end{equation}
	where $\mu$ is a spherical measure on $(G,K)$, is spherical for $(G,K)$. Then it is automatic that $\phi$ is bi-K-invariant and that $\phi(1)=1$.
	\end{definition}
	\begin{rem}
\normalfont
	By a representation $F:\mathit{L}^{1}\left(K\setminus G/K\right)\to\mathbb{C}$ of $G$ on a Banach space $\mathcal{H} =\mathit{L}^{2}\left(K\setminus G/K\right)$ we define the homomorphism $$f\mapsto F_{g}(f)=\int_{G}f(g)\phi(g^{-1})d\mu_{G}(g)$$ of $G$ into the group of non-singular bounded operators $\mathcal{B(H)}$, with the requirement that, for every $f\in\mathcal{H}$, the mapping $$f\mapsto F_{g}(f)$$ of $G$ into $\mathcal{H}$ is (strongly) continuous. For every compact set $K\subset G$ and every $f\in\mathcal{H}$ the function $F_{g}(f), \ g\in K$, is also compact and bounded in $\mathcal{H}$ by (Banach-Steinhaus theorem) this shows that
	\begin{equation}\label{key}
	\phi(g)=\left\|F_{g}(f)\right\|
	\end{equation}
	is bounded on every compact set, where $\phi(g)$ is lower semi-continuous on $G$, and satisfies
	\begin{equation}
	\phi(g_{1}g_{2})\leq \phi(g_{1})\phi(g_{2})
	\end{equation}
	for every two $g_{1},g_{2}\in K$; such a function will be called a \textbf{semi-norm} on $G$.
	\end{rem}
	\begin{prop}
\normalfont
	Given a semi-norm $\phi$ and for every $f\in\mathit{L}^{1}\left(K\setminus G/K\right)$, with
	\begin{equation}\label{key}
	\left\|f\right\|_{\phi}=\int_{G}\left|f(g)\right|\phi(g)d\mu_{G}(g),
	\end{equation}
	 $\mathit{L}^{1}\left(K\setminus G/K\right)$ is a topological algebra.
	\end{prop}
	\paragraph*{Proof:}  For all  $f_{1}, f_{2} \in\mathit{L}^{1}\left(K\setminus G/K\right)$, (5.5) implies
	\begin{flalign*}
	\left\|F_{g}(f_{1}*f_{2})(\phi)\right\|&\leq\left\|F_{g}\right\|\left\|f_{1}*f_{2}(\phi)\right\|\\
	&=\phi(g)\int_{G}\left|f_{1}*f_{2}(g)\right|d\mu(g)\\
	&=\phi(g)\int_{G}\int_{G}\left|f_{1}(g)f_{2}(y^{-1}g)\right|d\mu(g)d\mu(y)\\
	&\leq\phi(g)\int_{G}\int_{G}\left|f_{1}(g)\right|\left|f_{2}(y^{-1}g)\right|d\mu(g)d\mu(y)\\
	&=\phi(g)\left\|f_{1}(g)\right\|\int_{G}\left|f_{2}(y^{-1}g)\right|d\mu(y)\\
	&=\phi(g)\left\|f_{1}\right\|\left\|f_{2}\right\|.
	\end{flalign*}
	Hence,
	\begin{equation}\label{key}
	\left\|f_{1}*f_{2}\right\|_{\phi}\leq\left\|f_{1}\right\|_{\phi}\left\|f_{2}\right\|_{\phi}.
	\end{equation}
	To this end, the space $C_{c}\left(K\setminus G/K\right)$ of absolutely integrable functions with respect to $\phi(g)d\mu_{G}(g)$ can be considered in a natural way as a complete normed algebra under the convolution product. Of course, $\mathit{L}^{1}(G)$ is everywhere dense in $C_{c}\left(K\setminus G/K\right)$ by the definition of integrable functions.
	
	The norm $$\left\|f\right\|=\int_{G}\left|f(x)\right|dx$$ turns the group algebra into a normed vector space. Owing to the additional property $$\left\|f_{1}*f_{2}\right\|\leq\left\|f_{1}\right\|\left\|f_{2}\right\|,$$ The algebra $\mathit{C}_{c}^{\sharp}(G)$ is a closed subalgebra. $\ \ \ \ \square$
\section{Algebra of Lie Pseudogroups}
	Here, we shall follow \cite{Yumaguzhin} introduce Lie pseudogroups and their Lie algebras.
	\begin{definition}
\normalfont
	Let $M$ be a smooth manifold and let $\mathcal{G}$ be a collection of diffeomorphisms of open subsets of $M$ into $M$. $\mathcal{G}$ is called a pseudogroup if the following hold:

	\begin{enumerate}
		\item [(1)] $\mathcal{G}$ is closed under restriction: if $f\in\mathcal{G}$ and $U$ is a domain of $f$, then $f|_{V}\in\mathcal{G}$ for any open $V \subset U$.
		\item  [(2)] if $f : U \to M$ is a diffeomorphism, $U =\underset{\alpha}{\cup}U_{\alpha}$, and $f|_{U_{\alpha}}\in\mathcal{G}$, then $f\in\mathcal{G}$.
		\item [(3)] $\mathcal{G}$ is closed under inverse: if $f\in\mathcal{G}$, then $f^{-1}\in\mathcal{G}$.
		\item [(4)] $\mathcal{G}$ is closed under composition: $f:U\to M$ and $g :f(U) \to M$ both	belong to $\mathcal{G}$, then $g \circ f\in\mathcal{G}$.
		\item [(5)] The identity diffeomorphism $M \to M$ belongs to $\mathcal{G}$.

	\end{enumerate}	
	\end{definition}
	By $\mathcal{J}^{n}(M)$ we denote the manifold of all $n-$jets of all diffeomorphisms	of open subsets of $M$ into $M$. By $\mathcal{J}^{n}\mathcal{G}$ we denote the set of all $n-$jets of all diffeomorphisms belonging to $\mathcal{G}$.
	\begin{definition}
\normalfont
	A pseudogroup $\mathcal{G}$ is a Lie pseudogroup if there exists an integer $n\geq0$, called the order of $\mathcal{G}$, such that
	\begin{enumerate}
		\item [(1)] The set $\mathcal{J}^{n}\mathcal{G}$ is a smooth submanifold of $\mathcal{J}^{n}(M)$.
		\item [(2)] A diffeomorphism $f : U \to M$ belongs to $\mathcal{G}$ if and only if $\left[f\right]_{p}^{n}\in\mathcal{J}^{n}\mathcal{G}$ for all $p \in U$.

	\end{enumerate}
	\end{definition}

	The submanifold $\mathcal{J}^{n}\mathcal{G}$ of a Lie pseudogroup $\mathcal{G}$ is called a system of Partial Differential Equations defining $\mathcal{G}$.

	A pseudogroup $\mathcal{G}$ is transitive if for any $p_{1}$, $p_{2}$ in $M$ there exists $f \in \mathcal{G}$ such that $f(p_{1}) = p_{2}$.
	\subsection{Lie Algebra of Pseudogroups}
	Let $\mathcal{G}$ be a Lie pseudogroup acting on manifold $M$, $\beta$  a vector field
in $M$, and let $\varphi_{t}$ be the flow of $\beta$.
The vector field $\beta$ is $\mathcal{G}-$vector field if its flow consists of diffeomorphisms belonging to $\mathcal{G}$, that is $\beta\in\mathcal{G}$ for all $t$.

	\begin{prop}
\normalfont
		The set of all $\mathcal{G}-$vector fields is a Lie subalgebra in the Lie algebra of all vector fields in $M$.
	\end{prop}
	 Thus, the following hold:
	\begin{enumerate}
		\item [(1)] Suppose $\beta_{1}$ and $\beta_{2}$ are vector fields in $M$. Then $$\left[\beta_{1}^{n},\beta_{2}^{n}\right]=\left[\beta_{1},\beta_{2}\right]^{n}.$$
		\item [(2)] A vector field $\beta_{1}$ in $M$ is a $\mathcal{G}-$vector field if and only if the vector field $\beta^{n}$ is tangent to the equation $\mathcal{J}^{n}\mathcal{G}$.
	\end{enumerate}
The statement of the proposition follows from these facts.
	
The Lie algebra of all $\mathcal{G}-$vector fields is called the Lie algebra of $\mathcal{G}$. We	denote it by $\mathfrak{G}$.
\paragraph*{Proof:} The proof can be would have been completed if we can establish the following lemmas.
	\begin{lem}
\normalfont
	Suppose $\beta_{1}$ and $\beta_{2}$ are vector fields in $M$. Then $$\left[\beta_{1}^{n},\beta_{2}^{n}\right]=\left[\beta_{1},\beta_{2}\right]^{n}.$$
	\end{lem}

	\paragraph*{Proof.}
	Let $\varphi_{1t}$ be the flow of the vector field $\beta_{1}$ in $M$ and let $\varphi_{2t}$ be the flow of the vector field $\beta_{2}$ in $M$. Then the flow $\varphi_{t}^{n}$ is defined in $\mathcal{J}^{n}(M)$ by the formula  $$\beta_{1}^{n}=\varphi_{1t}\left(\left[f\right]_{p_{1}}^{n}\right)=\left[\varphi_{1t}\right]_{f(p_{1})}^{n}\cdot\left[f\right]_{p_{1}}^{n}=\left[\varphi_{1t}\circ f\right]_{p_{1}}^{n}$$ and $$\beta_{2}^{n}=\varphi_{2t}\left(\left[f\right]_{p_{2}}^{n}\right)=\left[\varphi_{2t}\right]_{f(p_{2})}^{n}\cdot\left[f\right]_{p_{2}}^{n}=\left[\varphi_{2t}\circ f\right]_{p_{2}}^{n},$$ then by induction when $n=1$, we have
	\begin{flalign*}
	 \left[\beta_{1},\beta_{2}\right]&=\left[\varphi_{1t}\left(\left[f\right]_{p_{1}}\right),\varphi_{2t}\left(\left[f\right]_{p_{2}}\right)\right]\\
	 &=\left[\left[\varphi_{1t}\right]_{f(p_{1})}\cdot\left[f\right]_{p_{1}},\left[\varphi_{2t}\right]_{f(p_{2})}\cdot\left[f\right]_{p_{2}}\right]\\
	&=\left[\left[\varphi_{1t}\circ f\right]_{p_{1}},\left[\varphi_{2t}\circ f\right]_{p_{2}}\right]\\
	&=\left[\varphi_{1t}\circ f\right]_{p_{1}}\left[\varphi_{2t}\circ f\right]_{p_{2}}-\left[\varphi_{2t}\circ f\right]_{p_{2}}\left[\varphi_{1t}\circ f\right]_{p_{1}}\\
	&=\beta_{1}\beta_{2}-\beta_{2}\beta_{1}.
	\end{flalign*}
	When $n=2$, we have
	\begin{flalign*}
	 \left[\beta_{1}^{2},\beta_{2}^{2}\right]&=\left[\varphi_{1t}\left(\left[f\right]_{p_{1}}^{2}\right),\varphi_{2t}\left(\left[f\right]_{p_{2}}^{2}\right)\right]\\
	 &=\left[\left[\varphi_{1t}\right]_{f(p_{1})}^{2}\cdot\left[f\right]_{p_{1}}^{2},\left[\varphi_{2t}\right]_{f(p_{2})}^{2}\cdot\left[f\right]_{p_{2}}^{2}\right]\\
	&=\left[\left[\varphi_{1t}\circ f\right]_{p_{1}}^{2},\left[\varphi_{2t}\circ f\right]_{p_{2}}^{2}\right]\\
	&=\left[\varphi_{1t}\circ f\right]_{p_{1}}^{2}\left[\varphi_{2t}\circ f\right]_{p_{2}}^{2}-\left[\varphi_{2t}\circ f\right]_{p_{2}}^{2}\left[\varphi_{1t}\circ f\right]_{p_{1}}^{2}\\
	&=\left(\left[\varphi_{1t}\circ f\right]_{p_{1}}\left[\varphi_{2t}\circ f\right]_{p_{2}}-\left[\varphi_{2t}\circ f\right]_{p_{2}}\left[\varphi_{1t}\circ f\right]_{p_{1}}\right)^{2}\\
	&=\left(\beta_{1}\beta_{2}-\beta_{2}\beta_{1}\right)^{2}\\
	&=\left[\beta_{1},\beta_{2}\right]^{2}.
	\end{flalign*}
	Now, when $n=n+1$ we have
	\begin{flalign*}
	 \left[\beta_{1}^{n+1},\beta_{2}^{n+1}\right]&=\left[\varphi_{1t}\left(\left[f\right]_{p_{1}}^{n+1}\right),\varphi_{2t}\left(\left[f\right]_{p_{2}}^{n+1}\right)\right]\\
	 &=\left[\left[\varphi_{1t}\right]_{f(p_{1})}^{n+1}\cdot\left[f\right]_{p_{1}}^{n+1},\left[\varphi_{2t}\right]_{f(p_{2})}^{n+1}\cdot\left[f\right]_{p_{2}}^{n+1}\right]\\
	&=\left[\left[\varphi_{1t}\circ f\right]_{p_{1}}^{n+1},\left[\varphi_{2t}\circ f\right]_{p_{2}}^{n+1}\right]\\
	&=\left[\varphi_{1t}\circ f\right]_{p_{1}}^{n+1}\left[\varphi_{2t}\circ f\right]_{p_{2}}^{n+1}-\left[\varphi_{2t}\circ f\right]_{p_{2}}^{n+1}\left[\varphi_{1t}\circ f\right]_{p_{1}}^{n+1}\\
	&=\left(\left[\varphi_{1t}\circ f\right]_{p_{1}}\left[\varphi_{2t}\circ f\right]_{p_{2}}\right)^{n+1}-\left(\left[\varphi_{2t}\circ f\right]_{p_{2}}\left[\varphi_{1t}\circ f\right]_{p_{1}}\right)^{n+1}\\
	&=\left(\beta_{1}\beta_{2}\right)^{n+1}-\left(\beta_{2}\beta_{1}\right)^{n+1}\\
	&=\left[\beta_{1},\beta_{2}\right]^{n+1}.
	\end{flalign*}
	Since this is true for $n=1$, $n=2$ and $n=n+1$, it is therefore true for all  $n$. \ \ \ \ \ $\square$
	
	The antisymmetric property of the Lie algebra follows directly from the following:
	\begin{lem}
\normalfont
		Suppose $\beta,\beta_{1},\beta_{2}$ are vector fields in $M$, for all $p,p_{1},p_{2}\in M$, then
		\begin{enumerate}
			\item [(1)] $\left[\beta,\beta\right]=0$ and
			\item [(2)] $\left[\beta_{1},\beta_{2}\right]=-\left[\beta_{2},\beta_{1}\right]$.
		\end{enumerate}
	\end{lem}
	\paragraph*{Proof.}
	(1). Let $\varphi_{t}$ be the flow of $\beta^{n}$ in $M$ with $\varphi_{t}$ defined in $\mathcal{J}^{n}(M)$ as $$\varphi_{t}\left(\left[f\right]_{p}^{n}\right)=\left[\varphi_{t}\right]_{f(p)}^{n}\cdot\left[f\right]_{p}^{n}=\left[\varphi_{t}\circ f\right]_{p}^{n}$$ for all $p\in M$. Then $\left[\beta^{n},\beta^{n}\right]$ implies
	\begin{flalign*}
	\left[\varphi_{t}\left(\left[f\right]_{p}^{n}\right),\varphi_{t}\left(\left[f\right]_{p}^{n}\right)\right]&=\left[\left[\varphi_{t}\circ f\right]_{p}^{n},\left[\varphi_{t}\circ f\right]_{p}^{n}\right]\\
	&=\left[\varphi_{t}\circ f\right]_{p}^{n}\left[\varphi_{t}\circ f\right]_{p}^{n}-\left[\varphi_{t}\circ f\right]_{p}^{n}\left[\varphi_{t}\circ f\right]_{p}^{n}\\
	 &=\left[\varphi_{t}\right]_{f(p)}^{n}\cdot\left[f\right]_{p}^{n}\left[\varphi_{t}\right]_{f(p)}^{n}\cdot\left[f\right]_{p}^{n}-\left[\varphi_{t}\right]_{f(p)}^{n}\cdot\left[f\right]_{p}^{n}\left[\varphi_{t}\right]_{f(p)}^{n}\cdot\left[f\right]_{p}^{n}\\
	&=\left[\varphi_{t}\right]_{f(p)}^{2n}\left[f\right]_{p}^{2n}-\left[\varphi_{t}\right]_{f(p)}^{2n}\left[f\right]_{p}^{2n}\\
	&=\left[\varphi_{t}\right]_{f(p)}^{2n}\left(\left[f\right]_{p}^{2n}-\left[f\right]_{p}^{2n}\right)=\left[\varphi_{t}\right]_{f(p)}^{2n}\cdot 0=0.
	\end{flalign*}
	(2).  We have
	\begin{flalign*}
	&\left[\varphi_{1t}\left(\left[f\right]_{p_{1}}^{n}\right),\varphi_{2t}\left(\left[f\right]_{p_{2}}^{n}\right)\right]\\
	&=\left[\left[\varphi_{1t}\circ f\right]_{p_{1}}^{n},\left[\varphi_{2t}\circ f\right]_{p_{2}}^{n}\right]\\
	&=\left[\varphi_{1t}\circ f\right]_{p_{1}}^{n}\left[\varphi_{2t}\circ f\right]_{p_{2}}^{n}-\left[\varphi_{2t}\circ f\right]_{p_{2}}^{n}\left[\varphi_{1t}\circ f\right]_{p_{1}}^{n}\\
	&=\left(\left[\varphi_{1t}\circ f\right]_{p_{1}}\left[\varphi_{2t}\circ f\right]_{p_{2}}-\left[\varphi_{2t}\circ f\right]_{p_{2}}\left[\varphi_{1t}\circ f\right]_{p_{1}}\right)^{n}\\
	 &=\left[\varphi_{1t}\right]_{f(p_{1})}^{n}\left[f\right]_{p_{1}}^{n}\left[\varphi_{2t}\right]_{f(p_{2})}^{n}\left[f\right]_{p_{2}}^{n}-\left[\varphi_{2t}\right]_{f(p_{2})}^{n}\left[f\right]_{p_{2}}^{n}\left[\varphi_{1t}\right]_{f(p_{1})}^{n}\left[f\right]_{p_{1}}^{n}\left[\varphi_{2t}\right]_{f(p_{2})}^{n}\left[f\right]_{p_{2}}^{n}\\
	 &=-\left[\varphi_{2t}\right]_{f(p_{2})}^{n}\left[f\right]_{p_{2}}^{n}\left[\varphi_{1t}\right]_{f(p_{1})}^{n}\left[f\right]_{p_{1}}^{n}\left[\varphi_{2t}\right]_{f(p_{2})}^{n}\left[f\right]_{p_{2}}^{n}+\left[\varphi_{1t}\right]_{f(p_{1})}^{n}\left[f\right]_{p_{1}}^{n}\left[\varphi_{2t}\right]_{f(p_{2})}^{n}\left[f\right]_{p_{2}}^{n}\\
	 &=-\left(\left[\varphi_{2t}\right]_{f(p_{2})}^{n}\left[f\right]_{p_{2}}^{n}\left[\varphi_{1t}\right]_{f(p_{1})}^{n}\left[f\right]_{p_{1}}^{n}\left[\varphi_{2t}\right]_{f(p_{2})}^{n}\left[f\right]_{p_{2}}^{n}-\left[\varphi_{1t}\right]_{f(p_{1})}^{n}\left[f\right]_{p_{1}}^{n}\left[\varphi_{2t}\right]_{f(p_{2})}^{n}\left[f\right]_{p_{2}}^{n}\right)\\
	&=-\left(\left[\varphi_{2t}\circ f\right]_{p_{2}}^{n}\left[\varphi_{1t}\circ f\right]_{p_{1}}^{n}-\left[\varphi_{1t}\circ f\right]_{p_{1}}^{n}\left[\varphi_{2t}\circ f\right]_{p_{2}}^{n}\right)\\
	 &=-\left(\varphi_{2t}\left(\left[f\right]_{p_{2}}^{n}\right)\varphi_{1t}\left(\left[f\right]_{p_{1}}^{n}\right)-\varphi_{1t}\left(\left[f\right]_{p_{1}}^{n}\right)\varphi_{2t}\left(\left[f\right]_{p_{2}}^{n}\right)\right)\\
	&=-\left[\varphi_{2t}\left(\left[f\right]_{p_{2}}^{n}\right),\varphi_{1t}\left(\left[f\right]_{p_{1}}^{n}\right)\right].
	\end{flalign*}
	Hence, $\left[\beta_{1}^{n},\beta_{2}^{n}\right]=-\left[\beta_{2}^{n},\beta_{1}^{n}\right]$. \ \ \ \ \ $\square$
	
	Next, let us show that the Jacobi identity also holds.
	\begin{lem}
\normalfont
	Suppose $\beta_{1}$, $\beta_{2}$ and $\beta_{3}$ are vector fields in $M$. Then $$\left[\beta_{1}^{n},\left[[\beta_{2}^{n},\beta_{3}\right]\right]+\left[\beta_{2}^{n},\left[\beta_{3}^{n},\beta_{1}^{n}\right]\right]+\left[\beta_{3}^{n},\left[\beta_{1}^{n},\beta_{2}^{n}\right]\right]=0.$$
	\end{lem}
	\paragraph*{Proof.}
	Let $\varphi_{1t}$ be the flow of $\beta_{1}$ in $M$, $\varphi_{2t}$ the flow of $\beta_{2}$ in $M$ and $\varphi_{3t}$ the flow of $\beta_{3}$ in $M$, where the flows $\varphi_{1t}$, $\varphi_{2t}$ and $\varphi_{3t}$ are defined in $\mathcal{J}^{n}(M)$ respectively by the formulae
	 $$\beta_{1}^{n}=\varphi_{1t}\left(\left[f\right]_{p_{1}}^{n}\right)=\left[\varphi_{1t}\right]_{f(p_{1})}^{n}\cdot\left[f\right]_{p_{1}}^{n}=\left[\varphi_{1t}\circ f\right]_{p_{1}}^{n}$$ and $$\beta_{2}^{n}=\varphi_{2t}\left(\left[f\right]_{p_{2}}^{n}\right)=\left[\varphi_{2t}\right]_{f(p_{2})}^{n}\cdot\left[f\right]_{p_{2}}^{n}=\left[\varphi_{2t}\circ f\right]_{p_{2}}^{n}$$ and
	 $$\beta_{3}^{n}=\varphi_{3t}\left(\left[f\right]_{p_{3}}^{n}\right)=\left[\varphi_{3t}\right]_{f(p_{3})}^{n}\cdot\left[f\right]_{p_{3}}^{n}=\left[\varphi_{3t}\circ f\right]_{p_{3}}^{n}$$ for all $p_{1}, p_{2}, p_{3}\in M$. Then $\left[\beta_{1},\left[\beta_{2},\beta_{3}\right]\right]$ implies
	\begin{flalign*}
	&\left[\left[\varphi_{1t}\circ f\right]_{p_{1}}^{n},\left[\left[\varphi_{2t}\circ f\right]_{p_{2}}^{n},\left[\varphi_{3t}\circ f\right]_{p_{3}}^{n}\right]\right]\\
	&=\left[\left[\varphi_{1t}\circ f\right]_{p_{1}}^{n},\left[\varphi_{2t}\circ f\right]_{p_{2}}^{n}\left[\varphi_{3t}\circ f\right]_{p_{3}}^{n}-\left[\varphi_{3t}\circ f\right]_{p_{3}}^{n}\left[\varphi_{2t}\circ f\right]_{p_{2}}^{n}\right]\\
	&=\left[\varphi_{1t}\circ f\right]_{p_{1}}^{n}\left[\varphi_{2t}\circ f\right]_{p_{2}}^{n}\left[\varphi_{3t}\circ f\right]_{p_{3}}^{n}-\left[\varphi_{3t}\circ f\right]_{p_{3}}^{n}\left[\varphi_{2t}\circ f\right]_{p_{2}}^{n}\left[\varphi_{1t}\circ f\right]_{p_{1}}^{n}\\
	&=\left[\varphi_{1t}\right]_{f(p_{1})}^{n}\cdot \left[f\right]_{p_{1}}^{n}\left[\varphi_{2t}\right]_{f(p_{2})}^{n}\cdot\left[f\right]_{p_{2}}^{n}\left[\varphi_{3t}\right]_{f(p_{3})}^{n}\cdot\left[f\right]_{p_{3}}^{n}\\
	&-\left[\varphi_{3t}\right]_{f(p_{3})}^{n}\cdot \left[f\right]_{p_{3}}^{n}\left[\varphi_{2t}\right]_{f(p_{2})}^{n}\cdot \left[f\right]_{p_{2}}^{n}\left[\varphi_{1t}\right]_{f(p_{1})}^{n}\cdot\left[f\right]_{p_{1}}^{n}\\
	 &=\left[\varphi_{1t}\right]_{f(p_{1})}^{n}\left[\varphi_{2t}\right]_{f(p_{2})}^{n}\left[\varphi_{3t}\right]_{f(p_{3})}^{n}\cdot\left(\left[f\right]_{p_{1}}^{n}\left[f\right]_{p_{2}}^{n}\left[f\right]_{p_{3}}^{n}\right)\\
	 &-\left[\varphi_{1t}\right]_{f(p_{1})}^{n}\left[\varphi_{2t}\right]_{f(p_{2})}^{n}\left[\varphi_{3t}\right]_{f(p_{3})}^{n}\cdot\left(\left[f\right]_{p_{1}}^{n}\left[f\right]_{p_{2}}^{n}\left[f\right]_{p_{3}}^{n}\right)\\
	 &=\left[\varphi_{1t}\right]_{f(p_{1})}^{n}\left[\varphi_{2t}\right]_{f(p_{2})}^{n}\left[\varphi_{3t}\right]_{f(p_{3})}^{n}\cdot\left(\left[f\right]_{p_{1}}^{n}\left[f\right]_{p_{2}}^{n}\left[f\right]_{p_{3}}^{n}-\left[f\right]_{p_{1}}^{n}\left[f\right]_{p_{2}}^{n}\left[f\right]_{p_{3}}^{n}\right)\\
	&=\left[\varphi_{1t}\right]_{f(p_{1})}^{n}\left[\varphi_{2t}\right]_{f(p_{2})}^{n}\left[\varphi_{3t}\right]_{f(p_{3})}^{n}\cdot 0=0.
	\end{flalign*}
	At the same time, $\left[\beta_{2},\left[\beta_{3},\beta_{1}\right]\right]$ implies
	\begin{flalign*}
	&\left[\left[\varphi_{2t}\circ f\right]_{p_{2}}^{n},\left[\left[\varphi_{3t}\circ f\right]_{p_{3}}^{n},\left[\varphi_{1t}\circ f\right]_{p_{1}}^{n}\right]\right]\\
	&=\left[\left[\varphi_{2t}\circ f\right]_{p_{2}}^{n},\left[\varphi_{3t}\circ f\right]_{p_{3}}^{n}\left[\varphi_{1t}\circ f\right]_{p_{1}}^{n}-\left[\varphi_{1t}\circ f\right]_{p_{1}}^{n}\left[\varphi_{3t}\circ f\right]_{p_{3}}^{n}\right]\\
	&=\left[\varphi_{2t}\circ f\right]_{p_{2}}^{n}\left[\varphi_{3t}\circ f\right]_{p_{3}}^{n}\left[\varphi_{1t}\circ f\right]_{p_{1}}^{n}-\left[\varphi_{1t}\circ f\right]_{p_{1}}^{n}\left[\varphi_{3t}\circ f\right]_{p_{3}}^{n}\left[\varphi_{2t}\circ f\right]_{p_{2}}^{n}\\
	 &=\left[\varphi_{2t}\right]_{f(p_{2})}^{n}\left[f\right]_{p_{2}}^{n}\left[\varphi_{3t}\right]_{f(p_{3})}^{n}\left[f\right]_{p_{3}}^{n}\left[\varphi_{1t}\right]_{f(p_{1})}^{n}\left[f\right]_{p_{1}}^{n}\\
	 &-\left[\varphi_{1t}\right]_{f(p_{1})}^{n}\left[f\right]_{p_{1}}^{n}\left[\varphi_{3t}\right]_{f(p_{3})}^{n}\left[f\right]_{p_{3}}^{n}\left[\varphi_{2t}\right]_{f(p_{2})}^{n}\left[f\right]_{p_{2}}^{n}\\
	 &=\left[\varphi_{1t}\right]_{f(p_{1})}^{n}\left[\varphi_{2t}\right]_{f(p_{2})}^{n}\left[\varphi_{3t}\right]_{f(p_{3})}^{n}\cdot\left(\left[f\right]_{p_{1}}^{n}\left[f\right]_{p_{2}}^{n}\left[f\right]_{p_{3}}^{n}\right)\\
	 &-\left[\varphi_{1t}\right]_{f(p_{1})}^{n}\left[\varphi_{2t}\right]_{f(p_{2})}^{n}\left[\varphi_{3t}\right]_{f(p_{3})}^{n}\cdot\left(\left[f\right]_{p_{1}}^{n}\left[f\right]_{p_{2}}^{n}\left[f\right]_{p_{3}}^{n}\right)\\
	 &=\left[\varphi_{1t}\right]_{f(p_{1})}^{n}\left[\varphi_{2t}\right]_{f(p_{2})}^{n}\left[\varphi_{3t}\right]_{f(p_{3})}^{n}\cdot\left(\left[f\right]_{p_{1}}^{n}\left[f\right]_{p_{2}}^{n}\left[f\right]_{p_{3}}^{n}-\left[f\right]_{p_{1}}^{n}\left[f\right]_{p_{2}}^{n}\left[f\right]_{p_{3}}^{n}\right)\\
	&=\left[\varphi_{1t}\right]_{f(p_{1})}^{n}\left[\varphi_{2t}\right]_{f(p_{2})}^{n}\left[\varphi_{3t}\right]_{f(p_{3})}^{n}\cdot 0=0.
	\end{flalign*}
	And lastly $\left[\beta_{3}^{n},\left[\beta_{1}^{n},\beta_{2}^{n}\right]\right]$ implies
	\begin{flalign*}
	&\left[\left[\varphi_{3t}\circ f\right]_{p_{3}}^{n},\left[\left[\varphi_{1t}\circ f\right]_{p_{1}}^{n},\left[\varphi_{2t}\circ f\right]_{p_{2}}^{n}\right]\right]\\
	&=\left[\left[\varphi_{3t}\circ f\right]_{p_{3}}^{n},\left[\varphi_{1t}\circ f\right]_{p_{1}}^{n}\left[\varphi_{2t}\circ f\right]_{p_{2}}^{n}-\left[\varphi_{2t}\circ f\right]_{p_{2}}^{n}\left[\varphi_{1t}\circ f\right]_{p_{1}}^{n}\right]\\
	&=\left[\varphi_{3t}\circ f\right]_{p_{3}}^{n}\left[\varphi_{1t}\circ f\right]_{p_{1}}^{n}\left[\varphi_{2t}\circ f\right]_{p_{2}}^{n}-\left[\varphi_{2t}\circ f\right]_{p_{2}}^{n}\left[\varphi_{1t}\circ f\right]_{p_{1}}^{n}\left[\varphi_{3t}\circ f\right]_{p_{3}}^{n}\\
	&=\left[\varphi_{3t}\right]_{f(p_{3})}^{n}\left[ f\right]_{p_{3}}^{n}\left[\varphi_{1t}\right]_{f(p_{1})}^{n}\left[ f\right]_{p_{1}}^{n}\left[\varphi_{2t}\right]_{f(p_{2})}^{n}\left[ f\right]_{p_{2}}^{n}\\
	&-\left[\varphi_{2t}\right]_{f(p_{2})}^{n}\left[ f\right]_{p_{2}}^{n}\left[\varphi_{1t}\right]_{f(p_{1})}^{n}\left[ f\right]_{p_{1}}^{n}\left[\varphi_{3t}\right]_{f(p_{3})}^{n}\left[ f\right]_{p_{3}}^{n}\\
	&=\left[\varphi_{3t}\right]_{f(p_{3})}^{n}\left[\varphi_{1t}\right]_{f(p_{1})}^{n}\left[\varphi_{2t}\right]_{f(p_{2})}^{n}\cdot\left(\left[ f\right]_{p_{3}}^{n}\left[ f\right]_{p_{1}}^{n}\left[ f\right]_{p_{2}}^{n}\right)\\
	 &-\left[\varphi_{2t}\right]_{f(p_{2})}^{n}\left[\varphi_{1t}\right]_{f(p_{1})}^{n}\left[\varphi_{3t}\right]_{f(p_{3})}^{n}\left(\left[f\right]_{p_{2}}^{n}\left[f\right]_{p_{1}}^{n}\left[f\right]_{p_{3}}^{n}\right)\\
	 &=\left[\varphi_{1t}\right]_{f(p_{1})}^{n}\left[\varphi_{2t}\right]_{f(p_{2})}^{n}\left[\varphi_{3t}\right]_{f(p_{3})}^{n}\cdot\left(\left[f\right]_{p_{1}}^{n}\left[f\right]_{p_{2}}^{n}\left[f\right]_{p_{3}}^{n}\right)\\
	 &-\left[\varphi_{1t}\right]_{f(p_{1})}^{n}\left[\varphi_{2t}\right]_{f(p_{2})}^{n}\left[\varphi_{3t}\right]_{f(p_{3})}^{n}\cdot\left(\left[f\right]_{p_{1}}^{n}\left[f\right]_{p_{2}}^{n}\left[f\right]_{p_{3}}^{n}\right)\\
	 &=\left[\varphi_{1t}\right]_{f(p_{1})}^{n}\left[\varphi_{2t}\right]_{f(p_{2})}^{n}\left[\varphi_{3t}\right]_{f(p_{3})}^{n}\cdot\left(\left[f\right]_{p_{1}}^{n}\left[f\right]_{p_{2}}^{n}\left[f\right]_{p_{3}}^{n}-\left[f\right]_{p_{1}}^{n}\left[f\right]_{p_{2}}^{n}\left[f\right]_{p_{3}}^{n}\right)\\
	&=\left[\varphi_{1t}\right]_{f(p_{1})}^{n}\left[\varphi_{2t}\right]_{f(p_{2})}^{n}\left[\varphi_{3t}\right]_{f(p_{3})}^{n}\cdot 0=0.
	\end{flalign*}
	Hence, $\left[\beta_{1},\left[\beta_{2},\beta_{3}\right]\right]+\left[\beta_{2},\left[\beta_{3},\beta_{1}\right]\right]+\left[\beta_{3},\left[\beta_{1},\beta_{2}\right]\right]=0. \ \ \ \ \ \square$\\
From the foregoing proposition 6.3 is established.
	\section{The Colombeau Algebra.}
	Distributions (see \cite{Colombeau},\cite{Oberguggenberger}, \cite{Koornwinder}) are presented by means of certain regularizations  with model delta nets, that is they will be considered as parametrized families $\left(f*\phi\right)_{\phi\in\mathcal{I}_{0}}$ where $\mathcal{I}_{0}$ is a subspace of $$\left\{\phi\in\mathcal{D}(\mathbb{R}^{n}):\int_{\mathbb{R}^{n}}\phi(t)dt=1\right\}.$$ The main focus naturally, will be on the subfamilies $$\left(f_{\epsilon}\right)_{\epsilon}=\left(f*\phi_{\epsilon}\right)_{\epsilon>0}$$ with $\phi_{\epsilon}$ given by $$\phi_{\epsilon}(t)=\epsilon^{-n}\phi\left(\frac{t}{\epsilon}\right).$$ To make the later differential-algebraic constructions work, we need to introduce an evaluation on the set of parameter.
	\begin{definition}
\normalfont
	Let $\phi\in\mathcal{D}(\mathbb{R}^{n})$ be a test function, we define the parameter $\mathcal{I}=\left(0,1\right)$ as follows:
	\begin{equation}
	\mathcal{I}(\mathbb{R}^{n})=\left\{\phi\in\mathcal{D}(\mathbb{R}^{n}):\int_{\mathbb{R}^{n}}\phi(t)dt=1\right\}
	\end{equation}
	and
	\begin{equation}
	\mathcal{I}(\mathbb{R}^{n})=\left\{\phi\in\mathcal{D}(\mathbb{R}^{n}):\int_{\mathbb{R}^{n}}t^{\alpha}\phi(t)dt=0\right\}
	\end{equation}
	for all  $|\alpha|\geq 1$,
	where $t=\left(t_{1},\dots,t_{n}\right)\in\mathbb{R}^{n}$, $\alpha=\left(\alpha_{1},\alpha_{2},\cdots,\alpha_{n}\right)\in\mathbf{N}^{n}$, and $t^{\alpha}=(t_{1})^{\alpha_{1}}\cdots(t_{n})^{\alpha_{n}}$.
	This parameter is called a mollifier.
	\end{definition}
	
	We have stated that for $\epsilon>0$, $\phi_{\epsilon}(t)=\frac{1}{\epsilon^{n}}\phi(\frac{1}{\epsilon})$. Now, if $\phi_{\epsilon}(t)\geq0$, then supp$\left(\phi_{\epsilon}(t)\right)=\overline{\mathbb{B}\left(0,\epsilon\right)}$, that is
	\begin{flalign*}
	\int_{\mathbb{R}^{n}}\phi_{\epsilon}(t)dt&=\int_{\mathbb{R}^{n}}\frac{1}{\epsilon^{n}}\phi(\frac{t}{\epsilon})dt\\
	&=\frac{1}{\epsilon^{n}}\int_{\mathbb{R}^{n}}\epsilon^{n}\phi(x)dx\\
	&=\int_{\mathbb{R}^{n}}\phi(x)dx=1=\int_{\mathbb{R}^{n}}\phi(t)dt.
	\end{flalign*}
	Here we have used the change of variable $t=\epsilon x$ and  $dt = \epsilon^{n}dx$. $\phi$.
	\begin{definition}
\normalfont
	We define \begin{flalign*}
	 \mathcal{E}\left(\mathbb{R}^{n}\right)&=\left(\mathcal{C}^{\infty}\left(\mathbb{R}^{n}\right)\right)^{\mathcal{I}\left(\mathbb{R}^{n}\right)}\\
	&=\left\{f:\mathcal{I}\left(\mathbb{R}^{n}\right)\to\mathcal{C}^{\infty}\left(\mathbb{R}^{n}\right)\right\}
	\end{flalign*}
	as the set of all $\mathcal{C}^{\infty}-$functions in $t$ for each fixed $\phi\in\mathcal{I}\left(\mathbb{R}^{n}\right)$.
	\end{definition}
	\begin{definition}
\normalfont
	The subset $\mathcal{E}_{M}\left(\mathbb{R}^{n}\right)$ of all $\left(f*\phi_{\epsilon}\right)_{\epsilon}=\left(f_{\epsilon}\right)_{\epsilon}\in\mathcal{E}\left(\mathbb{R}^{n}\right)$ such that: For all compact subsets $K$ of $\mathbb{R}^{n}$, for all $\alpha\in\mathbf{N}_{0}^{n}$ there exists $N\in\mathbb{N}$ such that the seminorm \begin{equation}\label{em}
	p_{N}\left(f\right)=\underset{t\in K}{\text{sup}}\left|\partial^{\alpha}f_{\epsilon}(t)\right|=\mathcal{O}\left(\epsilon^{-N}\right)\leq c\epsilon^{-N}
	\end{equation}
	as $\epsilon\to0$ holds, where $c>0$ and  $\phi\in\mathcal{I}\left(\mathbb{R}^{n}\right)$. The elements of $\mathcal{E}_{M}\left(\mathbb{R}^{n}\right)$ are called moderate. They also constitute a differential algebra.
	\end{definition}
	\begin{definition}
\normalfont
	The set $\mathcal{N}\left(\mathbb{R}^{n}\right)$ of all $\left(f_{\epsilon}\right)_{\epsilon}\in\mathcal{E}\left(\mathbb{R}^{n}\right)$ with the property that: For all compact subset $K$ of $\mathbb{R}^{n}$, $\alpha\in\mathbb{N}_{0}^{n}$ there exists $m\in\mathbb{N}$ such that the seminorm
	\begin{equation}\label{ne}
	p_{m}(f)=\underset{t\in K}{\text{sup}}\left|\partial^{\alpha}f_{\epsilon}(t)\right|=\mathcal{O}\left(\epsilon^{m}\right)\leq c\epsilon^{m}
	\end{equation}
	as $\epsilon\to0$ holds, where $c>0$ and  $\phi\in\mathcal{I}\left(\mathbb{R}^{n}\right)$. The elements of $\mathcal{N}\left(\mathbb{R}^{n}\right)$ are called neutral function and tends to zero faster than any power of $\epsilon$ when evaluated at $\phi_{\epsilon}$ with $\phi\in\mathcal{I}\left(\mathbb{R}^{n}\right)$ large enough. Clearly, $\mathcal{N}\left(\mathbb{R}^{n}\right)$ is a subalgebra closed under differentiation and it is an ideal of $\mathcal{E}_{M}\left(\mathbb{R}^{n}\right)$.
	\end{definition}
	\begin{definition}
\normalfont
	The algebra of generalized functions of Colombeau, denoted by $\mathcal{G}(\mathbb{R}^{n})$  (or $\mathcal{G}$),	 is the quotient algebra
	\begin{equation}\label{g}
	\mathcal{G}(\mathbb{R}^{n})=\mathcal{E}_{m}(\mathbb{R}^{n})/\mathcal{N}(\mathbb{R}^{n}).
	\end{equation}
	\end{definition}
	We remark that
	\begin{itemize}
		\item if $\overline{f}$ is a generalized function in $\mathcal{G}(\mathbb{R}^{n})$ then
		\begin{flalign*}
		\overline{f} &= \left(f*\phi_{\epsilon}\right)_{\epsilon} + \mathcal{N}(\mathbb{R}^{n})\\
		&=\left(f_{\epsilon}\right)_{\epsilon}+\mathcal{N}\left(\mathbb{R}^{n}\right)
		\end{flalign*}, where	 $\left(f_{\epsilon}\right)_{\epsilon}\in\mathcal{E}_{m}(\mathbb{R}^{n})$ is a representative of $\overline{f}$.
		\item Also, if $\overline{f} = \overline{g}$ in $\mathcal{G}(\mathbb{R}^{n})$ then
		\begin{equation*}
		 \left(f*\phi_{\epsilon}\right)_{\epsilon}-\left(g*\phi_{\epsilon}\right)_{\epsilon}=\left(f_{\epsilon}\right)_{\epsilon}-\left(g_{\epsilon}\right)_{\epsilon}\in\mathcal{N}(\mathbb{R}^{n})
		\end{equation*}
		where, $\left(f_{\epsilon}\right)_{\epsilon}$, $\left(g_{\epsilon}\right)_{\epsilon}$ are representatives of $\overline{f}$, $\overline{g}$ respectively.
	\end{itemize}
	\begin{lem}
\normalfont
	$\mathcal{G}(\mathbb{R}^{n})$ is an associative and a commutative algebra with identity.
	\end{lem}
	\paragraph*{Proof.} For all $\overline{f},\overline{g}\in\mathcal{G}\left(\mathbb{R}^{n}\right)$, we have
	\begin{flalign*}
	 \overline{f}\overline{g}(\phi)&=\left(\left(f*\phi_{\epsilon}\right)_{\epsilon}+\mathcal{N}\right)\cdot\left(\left(g*\phi_{\epsilon}\right)_{\epsilon}+\mathcal{N}\right)\\
	&=\left(f*\phi_{\epsilon}\right)_{\epsilon}\cdot\left(g*\phi_{\epsilon}\right)_{\epsilon}+\mathcal{N}\\
	&=\left(\left(fg\right)*\phi_{\epsilon}\right)_{\epsilon}+\mathcal{N}\\
	&=\left(f_{\epsilon}g_{\epsilon}\right)_{\epsilon}+\mathcal{N}\\
	&=\left(f_{\epsilon}\right)_{\epsilon}\left(g_{\epsilon}\right)_{\epsilon}+\mathcal{N}\\
	&\overline{fg}(\phi).
	\end{flalign*}
	It is obvious that $\partial^{\alpha}\mathcal{E}_{m}(\mathbb{R}^{n})\subset\mathcal{E}_{m}(\mathbb{R}^{n})$ and $\partial^{\alpha}\mathcal{N}(\mathbb{R}^{n})\subset\mathcal{N}(\mathbb{R}^{n})$, for all $\alpha$. Therefore, we can define
	\begin{equation}
	\partial^{\alpha}:\mathcal{G}(\mathbb{R}^{n})\to\mathcal{G}(\mathbb{R}^{n}):\overline{f}\mapsto\partial^{\alpha}\overline{f}
	\end{equation}
	where, \begin{flalign*}
	\partial^{\alpha}\overline{f}&=\partial^{\alpha}\left(f*\phi_{\epsilon}\right)+\mathcal{N}(\mathbb{R}^{n})\\
	&=\partial^{\alpha}f*\phi_{\epsilon}+\mathcal{N}(\mathbb{R}^{n}).
	\end{flalign*}
	It follows that $\partial^{\alpha}$ is linear, and satisfies Leibniz's rule of product derivatives. \ \ \ \ $\square$
	\begin{lem}
\normalfont
	$\mathcal{C}^{0}(\mathbb{R}^{n})$ is included in $\mathcal{E}_{M}\left(\mathbb{R}^{n}\right)$ as a linear subspace, not a subalgebra. Consequently, $\mathcal{C}^{0}\left(\mathbb{R}^{n}\right)$ is not a subalgebra of $\mathcal{G}\left(\mathbb{R}^{n}\right)$, either.
	\end{lem}
	\paragraph*{Proof.} In general, for all $f,g\in\mathcal{C}^{0}\left(\mathbb{R}^{n}\right)$ we have
	\begin{flalign*}
	 \left(f*\phi_{\epsilon}\right)\left(g*\phi_{\epsilon}\right)&=\int_{\mathbb{R}^{n}}f(t)\phi_{\epsilon}(t-x)dt\int_{\mathbb{R}^{n}}g(t)\phi_{\epsilon}(t-x)dt\\
	&=\frac{1}{\epsilon^{n}}\int_{\mathbb{R}^{n}}f(t)\phi(\frac{t-x}{\epsilon})dt\int_{\mathbb{R}^{n}}g(t)\phi(\frac{t-x}{\epsilon})dt\\
	&=\int_{\mathbb{R}^{n}}f(t+\epsilon y)\phi(y)dy\int_{\mathbb{R}^{n}}g(t+\epsilon y)\phi(y)dy\\
	&\neq\int_{\mathbb{R}^{n}}f\left(t+\epsilon y\right)g\left(t+\epsilon y\right)\phi(y)dy\\
	&=\int_{\mathbb{R}^{n}}fg(t+\epsilon y)\phi(y)dy\\
	&=fg*\phi_{\epsilon}.
	\end{flalign*}
	Another issue that arise here is that, if $f\in\mathcal{C}^{0}\left(\mathbb{R}^{n}\right)$, we can show that $\tilde{f_{1}}-\tilde{f_{2}}\in\mathcal{N}\left(\mathbb{R}^{n}\right)$, so both $\tilde{f_{1}}$ and $\tilde{f_{2}}$ are representatives of $f\in\mathcal{C}^{\infty}\left(\mathbb{R}^{n}\right)$. For convenience, we will show this in the case $n = 1$. Indeed, we have
	\begin{equation*}
	\left(\tilde{f_{1}}-\tilde{f_{2}}\right)\phi\left(t\right)=f(t)-\int_{\mathbb{R}}f\left(t+x\right)\phi(x)dx.
	\end{equation*}
	Therefore, one gets
	\begin{flalign*}
	\left(\tilde{f_{1}}-\tilde{f_{2}}\right)\left(\phi_{\epsilon}\right)&=f\left(t\right)-\int_{\mathbb{R}}f\left(t+\epsilon y\right)\phi\left(y\right)dy\\
	&=-\int_{\mathbb{R}}\left[f\left(t+\epsilon y\right)-f\left(t\right)\right]\phi\left(y\right)dy.
	\end{flalign*}
	Since $f\in\mathcal{C}^{\infty}\left(\mathbb{R}\right)$, we can apply Taylor's formula up to order $1$ to $f$ at the point $y$, and we get
	\begin{flalign*}
	f\left(t+\epsilon y\right)-f(t)&=(\epsilon y)\partial f"(t)+\frac{\left(\epsilon y\right)^{1+1}}{1+1}\partial f^{1+1}\left(t+\theta\epsilon y\right)\\
	&=\left(\epsilon y\right)f'(t)+\epsilon^{2}\frac{y^{2}}{2}\partial f^{(2)}\left(t+\theta\epsilon y\right)\\
	&=\left(\epsilon y\right)f'(t)+\epsilon^{2}\frac{y^{2}}{2}\partial f^{(2)}(t+\theta\epsilon y)
	\end{flalign*}
	where $0 < \theta < 1$. Hence, for arbitrary compact $K$ and $\phi\in\mathcal{I}$,	we have
	\begin{equation*}
	\left(\tilde{f_{1}}-\tilde{f_{2}}\right)\phi_{\epsilon}\left(t\right)=\mathcal{O}\left(\epsilon^{2}\right)
	\end{equation*}
	as $\epsilon\to 0$, uniformly on $K$ and $m=1$.
	
	Also, applying Taylor's formula up to order $m$ to $f$ at the point $y$, we obtain
	\begin{flalign*}
	f\left(t+\epsilon y\right)-f\left(t\right)&=\sum_{k=1}^{m}\frac{\left(\epsilon y\right)}{k!}\partial f^{\left(k\right)}(t)+\epsilon^{\left(m+1\right)}\frac{y^{m+1}}{\left(k+1\right)!}\partial f^{\left(k+1\right)}\left(t+\theta\epsilon y\right)\\
	&=\sum_{k=1}^{m}\frac{\left(\epsilon^{m} y\right)}{k!}\partial f^{(k)}(t)+\epsilon^{m+1}\frac{y^{m+1}}{\left(k+1\right)!}\partial f^{(k+1)}\left(t+\theta\epsilon y\right)
	\end{flalign*}
	where $0<\theta<1.$ Therefore for arbitrary $K\subset\mathbb{R}$ and for all $m\in\mathbb{N}$, we have
	$$\left(\tilde{f_{1}}-\tilde{f_{2}}\right)\phi_{\epsilon}(t)=\mathcal{O}\left(\epsilon^{m+1}\right)$$ as $\epsilon\to0$ uniformly on $K$ and $\phi\in\mathcal{I}$. This agrees with definition 7.4 for $|\alpha|\leq m$. This fact also holds for the estimate $\partial^{\alpha}\left(\tilde{f_{1}}-\tilde{f_{2}}\right)\phi_{\epsilon}\left(t\right)$, so $\tilde{f_{1}}-\tilde{f_{2}}\in\mathcal{N}\left(\mathbb{R}^{n}\right)$. \ \ \ \ $\square$
	\subsection{$\mathcal{L}^{1}\left(\mathbb{R}^{n}\right)$ Embedded in the	Colombeau Generalized Functions.}
	$\mathcal{L}^{1}\left(\mathbb{R}^{n}\right)\subset\mathcal{D}^{\prime}\left(\mathbb{R}^{n}\right)$ in the  distribution theory and hence, $\mathcal{L}^{1}\left(\mathbb{R}^{n}\right)-$functions are tempered distributions. So, to $f\in\mathcal{L}^{1}\left(\mathbb{R}^{n}\right)$, there corresponds an element in $\mathcal{G}\left(\mathbb{R}^{n}\right)$, denoted by $\tilde{f}+\mathcal{N}$ where $\tilde{f}\in\mathcal{E}_{M}\left(\mathbb{R}^{n}\right)$. Now, can we express $\tilde{f}$ in term of $f$ more clearly. In fact, we can do as follows:
	\begin{flalign*}
	\tilde{f}*\phi(t)&=\left\langle f(x),\phi(x-t)\right\rangle\\
	&=\int_{\mathbb{R}^{n}}f(x)\phi(x-t)dx,
	\end{flalign*}
	since $f\in\mathcal{L}^{1}\left(\mathbb{R}^{n}\right)$. So, we have $\tilde{f}+\mathcal{N}$ as the corresponding element of $f\in\mathcal{L}^{1}\left(\mathbb{R}^{n}\right)$. However, we cannot conclude from this that $\tilde{f}+\mathcal{N}_{\tau}\in\mathcal{G}_{\tau}\left(\mathbb{R}^{n}\right).$
	
	Since $f\in\mathcal{L}^{1}\left(\mathbb{R}^{n}\right)$, so the function $g$, where $$g(t) =\int_{-\infty}^{t}f(x)dx,$$ $t\in\mathbb{R}^{n}$ is continuous on $\mathbb{R}^{n}$. To this end,
	\begin{equation*}
	\left\|g(t)\right\|=\left|\int_{-\infty}^{t}f(y)dy\right|\leq\int_{\mathbb{R}^{n}}\left|f(y)\right|dy=\left\|f\right\|_{\mathcal{L}^{1}}
	\end{equation*}
	for all $t\in\mathbb{R}^{n}$.
	
	Hence, $g\in\mathcal{C}_{\tau}(\mathbb{R}^{n})$. In  $\mathcal{G}_{\tau}\left(\mathbb{R}^{n}\right)$, $g$ is assigned with $\tilde{g}+\mathcal{N}_{\tau}\left(\mathbb{R}^{n}\right)$, where $$\tilde{g}*\phi(t)=\int_{\mathbb{R}^{n}}g(t+x)\phi(x)dx.$$ It implies that $f$ belongs to $\mathcal{G}_{\tau}\left(\mathbb{R}^{n}\right)$, and it is assigned with $\partial\tilde{g}+\mathcal{N}_{\tau}\left(\mathbb{R}^{n}\right)$. Therefore, $f$ is assigned with the element $$\partial_{t}\left(\int_{\mathbb{R}^{n}}\left(\int_{-\infty}^{t+x}f(y)dy\right)\phi(x)dx\right)+\mathcal{N}_{\tau}\left(\mathbb{R}^{n}\right).$$
	
	Now, we notice that $\phi\in\mathcal{I}$ and $f\in\mathcal{L}^{1}\left(\mathbb{R}^{n}\right)$, we obtain
	\begin{flalign*}
	 \int_{-\infty}^{\infty}\left(\int_{-\infty}^{t+x}f(y)dy\right)\phi(x)dx&=\int_{-\infty}^{\infty}\left(\int_{-\infty}{t}f(t+x)dy\right)\phi(x)dx\\
	&=\int_{-\infty}^{t}\left(\int_{-\infty}^{\infty}\phi(x)f(y+x)dx\right)dy,
	\end{flalign*}
	and the inner integral as the function of $y$ is in $\mathcal{L}^{1}\left(\mathbb{R}^{n}\right)\bigcap\mathcal{C}^{\infty}\left(\mathbb{R}^{n}\right)$. it follows that
	\begin{flalign*}
	\int_{-\infty}^{\infty}\left(\int_{-\infty}^{t+x}f(y)dy\right)\phi(x)dx=\int_{-\infty}^{\infty}\phi(x)f(t+x)dx.
	\end{flalign*}
	Therefore, $f$ is assigned with the element $$\int_{\mathbb{R}^{n}}f(t + x)\phi(x)dx + \mathcal{N}_{\tau}$$ in $\mathcal{G}_{\tau}\left(\mathbb{R}^{n}\right)$.	 It also shows us that in $\mathcal{G}\left(\mathbb{R}^{n}\right)$, the function $f\in\mathcal{L}^{1}\left(\mathbb{R}^{n}\right)$ is assigned with the element $$\int_{\mathbb{R}^{n}}f(t+x)\phi(x)dx+\mathcal{N}.$$
	
	Now, we will use the results above to study the relationship between the	integral of $f \in \mathcal{L}^{1}\left(\mathbb{R}^{n}\right)$ in the usual sense and the one in the sense of tempered 	generalized function.
	\begin{lem}
\normalfont
	The topology of $\mathcal{G}\left(\mathbb{R}^{n}\right)$ is the topology it inherited as an embedding image of $\mathcal{L}^{1}\left(\mathbb{R}^{n}\right)$.
	\end{lem}
	\paragraph*{Proof.} It suffices to show that for all $\overline{f}\in\mathcal{G}\left(\mathbb{R}^{n}\right)$, such that $\overline{f}=\tilde{f}+\mathcal{N}\left(\mathbb{R}^{n}\right)$, where $\tilde{f}$ is the representative of $f$ in  $\mathcal{E}_{M}\left(\mathbb{R}^{n}\right)$, for all $\phi\in\mathcal{I}\left(\mathbb{R}^{n}\right)$, and $\epsilon$ small enough we have $$\left\|\tilde{f}*\phi_{\epsilon}(t)\right\|_{\mathcal{E}_{M}}\leq\mathit{c}\left\|f\right\|_{\mathcal{L}^{1}\left(\mathbb{R}^{n}\right)}.$$ To this end, we have
	\begin{flalign*}
	\left\|\tilde{f}*\phi_{\epsilon}(t)\right\|_{\mathcal{E}_{M}}&=\left|\int_{\mathbb{R}^{n}}f(t)\phi_{\epsilon}(t-x)dt\right|\\
	&\leq\int_{\mathbb{R}^{n}}\left|f(t)\frac{1}{\epsilon^{n}}\phi\left(\frac{t-x}{\epsilon}\right)\right|dt\\
	&\leq\epsilon^{-|n|}\int_{\mathbb{R}^{n}}\left|f(x+\epsilon y)\epsilon^{n}\phi(y)\right|dy\\
	&\leq\epsilon^{|-n|}\cdot\epsilon^{|n|}\int_{\mathbb{R}^{n}}\left|f(x+\epsilon y)\right|\left|\phi(y)\right|dy\\
	&\leq\mathit{c}\int_{\mathbb{R}^{n}}\left|f(x+\epsilon y)\right|dy\\
	&\leq\mathit{c}\left\|f\right\|_{\mathcal{L}^{1}}.
	\end{flalign*}
	Therefore, $$\left\|\tilde{f}*\phi_{\epsilon}(t)\right\|_{\mathcal{E}_{M}}=\mathcal{O}\left(\epsilon^{0}\right)\leq\mathit{c}\left\|f\right\|_{\mathcal{L}^{1}}$$ which means $\tilde{f}\in\mathcal{N}\left(\mathbb{R}^{n}\right)$ for all $n\in\mathbb{N}$.
	
	Also, if $$\widetilde{fg}+\mathcal{N}=\tilde{f}\tilde{g}+\mathcal{N}=\left(\tilde{f}+\mathcal{N}\right)\left(\tilde{g}+\mathcal{N}\right)$$ where $\tilde{f}$ and $\tilde{g}$ are representatives of $f$ and $g$ in $\mathcal{E}_{M}\left(\mathbb{R}^{n}\right)$ we have
	\begin{flalign*}
	 \left\|\left(\widetilde{fg}\right)*\phi_{\epsilon}(t)\right\|_{\mathcal{E}_{M}}&=\left\|\left(\tilde{f}\tilde{g}\right)*\phi_{\epsilon}(t)\right\|_{\mathcal{E}_{M}}\\
	&=\left\|\left(\tilde{f}*\phi_{\epsilon}(t)\right)\left(\tilde{g}*\phi_{\epsilon}(t)\right)\right\|_{\mathcal{E}_{M}}\\
	&\leq\left\|\tilde{f}*\phi_{\epsilon}(t)\right\|_{\mathcal{E}_{M}}\left\|\tilde{g}*\phi_{\epsilon}(t)\right\|_{\mathcal{E}_{M}}\\
	&\leq\mathit{c}_{1}\left\|f\right\|_{\mathcal{L}^{1}}\left\|\tilde{g}*\phi_{\epsilon}(t)\right\|_{\mathcal{E}_{M}}\\
	&\leq\mathit{c}_{1}\mathit{c}_{2}\left\|f\right\|_{\mathcal{L}^{1}}\left\|g\right\|_{\mathcal{L}^{1}}\\
	&\leq\mathit{C}\left\|f\right\|_{\mathcal{L}^{1}}\left\|g\right\|_{\mathcal{L}^{1}},
	\end{flalign*}
	where we have taking $\mathit{c}_{1}\mathit{c}_{2}=\mathit{C}$,
	for all $\phi\in\mathcal{I}\left(\mathbb{R}^{n}\right)$, $0<\epsilon<1$ and $t\in\mathbb{R}^{n}$.\ \ \ \ $\square$
\begin{center}

\end{center}
\end{document}